\documentclass{article}

\usepackage{latexsym,amssymb}
\title{\bf A generalized character related to the $p$-local structure and representation theory of a finite group}
\author{Geoffrey R. Robinson\\
Institute of Mathematics\\
University of Aberdeen\\
Fraser Noble Building\\
Aberdeen\\
AB24 3UE\\
United Kingdom\\}

\begin{document}

\maketitle

\medskip
\begin{center}
\emph{To the memory of Marty Isaacs}
\end{center}

\begin{abstract}

\medskip
Let $G$ be a finite group and $p$ be a prime.  Let $(\mathbb{K},R,\mathbb{F})$ be a (splitting) $p$-modular system for $G$. 

\medskip
We let $G_{p}$ denote the set of $p$-elements of $G$ and $G_{p^{\prime}}$ denote the set of $p$-regular elements of $G$. In this note, we examine the generalized character  $\Psi_{1,p,G}$ of $G$ which vanishes on all $p$-singular elements and whose value at each $p$-regular $y \in G$ is the number of $p$-elements of $C_{G}(y).$
We examine cases when $\Psi_{1,p,G}$ is a character, and when it is a character afforded by a projective $RG$-module. We conjecture that, in fact, $\Psi_{1,p,G}$ is always a character, and may always be afforded by a projective $RG$-module. We discuss some properties of $\Psi_{1,p,G}$, and relate it to a truncated version of the $RG$-module afforded by the conjugation action of $G$ on itself.

\medskip
Since the virtual character $\Psi_{1,p,G}$ vanishes on $p$-singular elements, it corresponds to a unique virtual projective module $P_{1,p,G}$ in the Green ring for $RG$, and we obtain a uniform  explicit expression for $P_{1,p,G}$ which is valid for all finite $G$. It is however not immediately obvious from this expression whether or not $P_{1,p,G}$  is  always a module. On the other hand, letting $R$ denote the trivial $RG$-module, we give a general explicit formula expressing $P_{1,p,G}-R$ (in the Green ring for $RG$) as a virtual module induced from (proper) $p$-local subgroups of $G$ when $O_{p}(G) = 1.$

\medskip
Among other examples, we prove for each prime $p$ and each prime power $q$ (not necessarily divisible by $p$) that $P_{1,p,G}$ is a genuine projective $RG$-module when $G \cong {\rm PSL}(2,q)$ and when $G \cong {\rm SL}(2,q).$ 
\end{abstract}

\section{Notation and Introduction} %1

\medskip
Let $G$ be a finite group, $p$ be a prime, $(\mathbb{K},R, \mathbb{F})$ be a (splitting) $p$-modular system for $G$. 

\medskip
When $n$ is a positive integer and $\pi$ is a set of primes, we let $n_{\pi}$ denote the largest positive integer which divides $n$ and is itself only divisible by primes in $\pi.$ As usual, when $g \in G$ and $\pi$ is a set of primes, we may (uniquely) write $g = xy = yx$ with $\langle x \rangle $ a $\pi$-group and $\langle y \rangle$ a $\pi^{\prime}$-group. We call $x$ the \emph{$\pi$-part} of $g$ and denote it by $g_{\pi}.$ The \emph{$\pi$-section of $g \in G$} is the set of elements of $G$ whose $\pi$-part is conjugate to $g_{\pi}$, and is denoted by ${\rm S}_{\pi}^{G}(g)$ (usually we will take $g$ to be a $\pi$-element itself).  The element $g \in G$ is said to be \emph{$\pi$-singular} if $g_{\pi} \neq 1_{G}$, and \emph{$\pi$-regular} otherwise.

\medskip
We let $G_{p}$ denote the set of $p$-elements of $G$ and $G_{p^{\prime}}$ denote the set of $p$-regular elements of $G$.  Let $\{ \phi_{i}: 1 \leq i \leq \ell \}$ be the irreducible Brauer characters of $G$, and let $\{ \theta_{i}: 1 \leq i \leq \ell \}$ be the set of characters of the projective indecomposable $RG$-modules, where we label so that $\langle \theta_{i}, \phi_{j} \rangle = \delta_{ij}$ for $1 \leq i,j \leq \ell.$ (When $\alpha$ is a class function which is only defined for $p$-regular elements and $\beta$ is a class function of $G$ which vanishes on $p$-singular elements, we will use $\langle \alpha, \beta \rangle$
to denote $$\frac{1}{|G|} \left( \sum_{g \in G_{p^{\prime}}} \alpha(g) \overline{\beta(g)} \right),$$ which is compatible with the usual notation for the standard inner product on class functions of $G$).

\medskip
We define the class function 
$\Psi_{1,p,G}$ of $G$ to take value $0$ on all $p$-singular elements of $G$, and to take value

\medskip
 $\Psi_{1,p,G}(y) =$ \emph{(the number of $p$-elements of $C_{G}(y)$)} whenever $y \in G$ is $p$-regular. 

\medskip
Notice that $\Psi_{1,p,G}$ agrees with the permutation character $$\sum_{ x \in G_{p}/G} {\rm Ind}_{C_{G}(x)}^{G}(1)$$ on $p$-regular elements (and vanishes elsewhere). We note also that it is clear that algebraically conjugate irreducible characters occur with equal multiplicity in $\Psi_{1,p,G}.$ 

\medskip
By a well-known theorem of Frobenius, we know that $|C_{G}(y)|_{p}$ divides $|C_{G}(y)_{p}|$ for each $p$-regular $y \in G$. By Brauer's characterization of characters, it follows that $\Psi_{1,p,G}$ is a generalized character of $G$.

\medskip
To see that $\Psi_{1,p,G}$ is a generalized character, notice that whenever $E$ is a Brauer elementary subgroup of $G$ of order prime to $p$,  the class function ${\rm Res}^{G}_{E}(\Psi_{1,p,G})$ is a character of $E$ (since it agrees with the restriction to $E$ of a permutation character). If a Brauer elementary subgroup $E$ has the form  $P \times F$ for some non-trivial $p$-group $P$ and some $p^{\prime}$-group $F$, then $${\rm Res}^{G}_{E}(\Psi_{1,p,G}) = {\rm Ind}_{F}^{E}( \frac{{\rm Res}^{G}_{F}(\Psi_{1,p,G})}{|P|})$$ and $$\frac{{\rm Res}^{G}_{F}( \Psi_{1,p,G})}{|P|}$$ is a character of $F.$ Using Brauer's induction theorem, it also follows that $\Psi_{1,p,G}$ is a $\mathbb{Z}$-linear combination of characters each induced from linear characters of Brauer elementary $p^{\prime}$-subgroups of $G.$

\medskip
Since $\Psi_{1,p,G}$ agrees with a character on $p$-regular elements, and each (irreducible character $\chi$ in a) $p$-block of defect zero of $G$ vanishes off $p$-regular elements, it follows that $\langle \Psi_{1,p,G}, \chi \rangle$ is a non-negative integer whenever $\chi$ is (an irreducible character in) a $p$-block of defect zero of $G$. In fact, we have 
$$\langle \Psi_{1,p,G}, \chi \rangle = \sum_{x \in G_{p}/G} \langle {\rm Res}^{G}_{C_{G}(x)}(\chi),1 \rangle$$ for each such $\chi.$

\medskip
Notice that if $x \in G_{p}$ and $z^{p^{k}} = x$, we have $C_{G}(z) \leq C_{G}(x)$, and that, furthermore,
two such $p^{k}$-th roots of $x$ which are conjugate in $G$ are already conjugate in $C_{G}(x).$

\medskip
Now if $\langle {\rm Res}_{C_{G}(x)}^{G}(\chi), 1 \rangle > 0$ we must have $\langle {\rm Res}_{C_{G}(z)}^{G}(\chi), 1 \rangle > 0$ for any such $z$. It follows by Frobenius reciprocity that $\langle \Psi_{1,p,G}, \chi \rangle$
is at least as great as the number of $C_{G}(x)$-conjugacy classes of $p$-power roots of $x$ when such an $x$ exists.

\medskip
In the case $p = 2,$ whenever $\chi$ is a real-valued irreducible character in a $2$-block of defect zero of $G$, there is (by Theorem 8 of Murray [12]) an involution\\ $t \in G$, unique up to conjugacy, such that 
$\langle \chi, {\rm Ind}_{C_{G}(t)}^{G}(1) \rangle > 0,$ and the multiplicity is in fact $1$. Hence the above argument shows that 
$\langle \Psi_{1,p,G}, \chi \rangle$ is at least as great as the number of $C_{G}(t)$-conjugacy classes of $2$-power
roots of $t$ (including $t$ itself).

\section{ The generalized character $\Psi_{1,p,G}$ and root counting}%2

\medskip
Notice that if $G$ has a Sylow $p$-subgroup $S$, then $\Psi_{1,p,G}(x)$ is the number of $|S|$-th roots of $x$ in $G$. If $x$ is $p$-singular, there is  clearly no such root of $x$, while if $x$ is $p$-regular and $y^{|S|} =x, $ then $y$ has the form $ uv = vu$ where $u$ is the unique generator of $\langle x \rangle$ with $u^{|S|} = x$ and $v$ is a $p$-element of $C_{G}(x).$ In fact, we may choose a power $q$ of $p$ so that $|S|$ divides $q$ and $q \equiv 1$ (mod $|G|_{p^{\prime}}$), and then we see easily that $\Psi_{1,p,G}(x)$ is the number of $q$-th roots of $x$ in $G$ for each $x \in G$.

\medskip
Now for each irreducible character $\chi$ of $G$, it is well-known that the class function $\chi^{(q)}$ of $G$ defined by $\chi^{(q)}(g) = \chi(g^{q})$
for all $g \in G$ is a generalized character of $G$.  Note that we have $\chi^{(q)}(g) = \chi(g_{p^{\prime}})$ for each $g$ by the choice of $q$.

\medskip
Setting $\nu_{q}(\chi) = \langle \chi^{(q)}, 1 \rangle,$ we see easily that $$\Psi_{1,p,G} = \sum_{\chi \in {\rm Irr}(G)}\nu_{q}(\chi) \chi,$$ so that $\Psi_{1,p,G}$ is a character of $G$ if and only if $\nu_{q}(\chi) \geq 0$ for each irreducible character $\chi$ of $G$.

\medskip
We may note further that $\nu_{q}(\chi) \leq \chi(1)$ with equality if and only if $$\chi(g_{p^{\prime}}) = \chi(1)$$ for all $g \in G,$ so that $\nu_{q}(\chi) = \chi(1)$ if and only if $O^{p}(G) \leq {\rm ker} \chi.$ 

\medskip
We also note that $$ -\chi(1) < \nu_{q}(\chi)  \leq \chi(1)$$ for each irreducible character $\chi$ of $G,$  so that when $\chi$ is linear, we have 
$\nu_{q}(\chi) \geq 0$ and $\nu_{q}(\chi) = 1 = \chi(1)$ if and only if $O^{p}(G) \leq {\rm ker} \chi.$ Hence we have proved:

\medskip
\noindent {\bf Theorem 2.1:} \emph{ We have $G \neq O^{p}(G)$ if and only if $\Psi_{1,p,G}$ contains some non-trivial linear character with non-zero  multiplicity. Furthermore, for every irreducible character $\chi$ of $G$, we have $$-\chi(1) < \nu_{q}(\chi) \leq \chi(1),$$ and the irreducible characters $\chi$ of $G$ which occur with multiplicity $\chi(1)$ in $\Psi_{1,p,G}$ are precisely the irreducible characters of $G$ with $O^{p}(G)$ in their kernels.}

\medskip
Another extreme case is: 

\medskip
\noindent {\bf Theorem 2.2:} \emph{ Let $\chi$ be an irreducible character of $G$ such that $$O^{p^{\prime}}(G) \leq {\rm ker} \chi.$$ Then $\chi$ occurs with zero multiplicity in $\Psi_{1,p,G}$ except when $\chi$ is the trivial character, in which case, $\chi$ occurs with multiplicity one in $\Psi_{1,p,G}.$ In any case, $\chi$ occurs with non-negative multiplicity in $\Psi_{1,p,G}.$}

\medskip
\noindent {\bf Proof:} The given hypothesis on $\chi$ is equivalent to asserting that every $p$-element of $G$ lies in ${\rm ker} \chi.$ Hence for each $g \in G,$ we have $\chi^{(q)}(g) = \chi(g_{p^{\prime}}) = \chi(g).$ Thus $\chi^{(q)} = \chi$, so that $\nu_{q}(\chi) = \delta_{1,\chi}$ and the result follows.

\medskip
\noindent {\bf Remark 2.3:} Note that Theorem 2.1 implies, in particular that $\Psi_{1,p,G}$ is the regular character of $G$ if and only if $G$ is a $p$-group. On the other hand, if $G$ is a $p^{\prime}$-group, then Theorem 2.2 implies that $\Psi_{1,p,G}$ is the trivial character. Finally, we note that if $\Psi_{1,p,G}$ is the trivial character, then $G$ contains only one $p$-element, so that (by Cauchy's theorem), $G$ is a $p^{\prime}$-group.

Hence we have:

\noindent {\bf Corollary 2.4:} \emph{ $\Psi_{1,p,G}$ is the trivial character if and only if $G$ is a $p^{\prime}$-group. Also, $\Psi_{1,p,G}$ is the regular character if and only if $G$ is a $p$-group.}

\medskip
\noindent {\bf Remark 2.5:} In the case of symmetric groups, generalized characters which count the numbers of roots of a given element are well-understood. T. Scharf has proved in [15] that for any two positive integers $m$ and $n$, the integer-valued class function $\theta_{m}$ of the symmetric group $S_{n}$ defined by letting $\theta_{m}(x)$ denote the number of $m$-th roots of $x$ in $S_{n}$ is a character of $S_{n}.$ No such result is true for the general finite group $G$. Notice that Scharf's result implies, in particular, that $\Psi_{1,p,S_{n}}$ is a character of $S_{n}$ for every prime $p$ and every positive integer $n$.

\section{Some general observations}%3

\medskip
\noindent {\bf Remark 3.1:} For a general finite group $G$, we note that $$\langle \Psi_{1,p,G},1 \rangle = \frac{1}{|G|} \left(\sum_{ y \in G_{p^{\prime}}/G} |S_{p^{\prime}}^{G}(y)|\right) = 1,$$ since for each $p$-regular $y \in G,$ the cardinality of the $p^{\prime}$-section of $y$ in $G$ is

\medskip
 $[G:C_{G}(y)]$ $\times$ (\emph{the number of $p$-elements of $C_{G}(y)$)}. 

\medskip
Hence, by block orthogonality relations, the character of the projective cover of the trivial module occurs exactly once when $\Psi_{1,p,G}$ is uniquely expressed as a $\mathbb{Z}$-linear combination of characters of projective indecomposable $RG$-modules.

\medskip
In this note, we are principally concerned with two questions: firstly, is $\Psi_{1,p,G}$ always a character? Secondly, is $\Psi_{1,p,G}$ a non-negative integer combination of characters of projective indecomposable $RG$-modules (which certainly would imply that it is a character)? We conjecture that the latter is always the case. Though we are not able to prove this at present, we illustrate that it is true in many cases, and we draw some consequences for groups such that $\Psi_{1,p,G}$ is a character afforded by a projective $RG$-module.

\medskip
We now note the following extension of Corollary 2.4: 

\medskip
\noindent {\bf Theorem 3.2 :} i)\emph{ Suppose that $G$ has a normal $p$-complement. Then $\Psi_{1,p,G}$ is a character, and may be afforded by a projective $RG$-module.}

\medskip
\noindent ii) \emph{ Suppose that $G$ has a normal Sylow $p$-subgroup. Then $\Psi_{1,p,G}$ is the character afforded by the projective cover of the trivial $RG$-module.}

\medskip
\noindent iii) \emph{ If $\Psi_{1,p,G/O_{p}(G)}$ is a character of $G/O_{p}(G)$ afforded by a projective $RG/O_{p}(G)$-module $M$, then $\Psi_{1,p,G}$ is afforded by $P(M),$ the projective cover (as $RG$-module) of the inflation of $M$ to an $RG$-module.}

\medskip
\noindent {\bf Proof:} i)  It suffices to prove that whenever $\phi_{i}$ is a Brauer irreducible character of $G$,
then $$\langle \Psi_{1,p,G}, \phi_{i} \rangle \geq 0.$$ But  $$\langle \Psi_{1,p,G}, \phi_{i} \rangle = \sum_{x \in G_{p}/G} \langle {\rm Res}^{G}_{C_{G}(x)}(\phi_{i}),1 \rangle^{\prime} \geq 0,$$ where $\langle,\rangle^{\prime}$ denotes that the inner product is restricted to $p$-regular elements.

\medskip
 For each $p$-element $x \in G$, the set of $p$-regular elements of $C_{G}(x)$ is the subgroup 
$O_{p^{\prime}}(C_{G}(x))$, and it follows that $$|C_{G}(x)|_{p}\langle {\rm Res}^{G}_{C_{G}(x)}(\phi_{i}),1 \rangle^{\prime}$$ is a non-negative integer, and, in particular, $\langle \Psi_{1,p,G}, \phi_{i} \rangle$ is non-negative, since $\phi_{i}$ always restricts to a character of $O_{p^{\prime}}(C_{G}(x))$ (and we already know that the inner product in question is an integer).

\medskip
To be precise, we have $$\Psi_{1,p,G} =  \sum_{i=1}^{\ell} \left(\sum_{ x \in G_{p}/G} \frac{\langle {\rm Res}^{G}_{O_{p^{\prime}}(C_{G}(x))}( \phi_{i}),1 \rangle}{|C_{G}(x)|_{p}}\right) \theta_{i} .$$ 

\medskip
\noindent ii) Suppose that $G$ has a normal Sylow $p$-subgroup $P$. Then $G$ has a Hall $p^{\prime}$-subgroup $H$, and the character ${\rm Ind}_{H}^{G}(1)$ is easily seen to take value $|C_{P}(y)|$ for each $y \in H$. But for each such $y$, the subgroup $C_{G}(y)$ has a normal Sylow $p$-subgroup $C_{P}(y)$ and we have $|C_{P}(y)| = |C_{G}(y)_{p}| = \Psi_{1,p,G}(y).$ Hence $$\Psi_{1,p,G} = {\rm Ind}_{H}^{G}(1),$$ which is the character afforded by the projective cover of the trivial module.

\medskip
\noindent iii) A theorem of W.F. Reynolds [13] asserts that if $H$ is a finite group and\\ $U = O_{p}(H)$, then whenever $\alpha$ is the Brauer character of a projective\\ $RH/U$-module $X$ and $\beta$ is the Brauer character afforded by the projective cover of $X$ as $RH$-module, then for each $p$-regular $y \in H,$ we have $$\beta(y) = |C_{U}(y)|\alpha(yU)$$ for all $p$-regular $y \in H$. In fact, Reynolds' observation gives a (Brauer character preserving) bijection between virtual Brauer characters of virtual projective $RH/U$-modules and virtual Brauer characters of virtual projective $RH$-modules, obtained by multiplying the virtual projective Brauer character of $H/U$ by the Brauer character afforded by the conjugation action of $H$ on $U$.

\medskip
\noindent {\bf Remark 3.3 :} Theorem 3.2 tells us that if $G$ has a normal $p$-complement or a normal Sylow $p$-subgroup, or if $\Psi_{1,p,G/O_{p}(G)}$ is a character afforded by a projective $RG$-module, then $\Psi_{1,p,G}$ is a character afforded by a projective $RG$-module. 
It follows that $\Psi_{1,p,G}$ is a character afforded by a projective $RG$-module if $G = O_{p,p^{\prime},p}(G).$ 

\medskip
By P. Fong's theorems on characters of $p$-solvable groups (see [7]), it is known that when $G$ is $p$-solvable, a character which vanishes on $p$-singular elements may be afforded by a projective $RG$-module. Hence to prove for such $G$ that $\Psi_{1,p,G}$ is a character which may be afforded by a projective $RG$-module, it is sufficient to prove that $\Psi_{1,p,G}$ is a character of $G$.

\medskip
\section{ The generalized character $\Psi_{1,p,G}$ and blocks of $RG$}%4
\medskip
Since $\Psi_{1,p,G}$ is a virtual projective character of $G$, it may be uniquely decomposed as an orthogonal sum $\sum_{b} \Psi^{(b)}_{1,p,G}$ where $b$ runs over  the blocks of $RG$, and each $\Psi^{(b)}_{1,p,G}$ is a (possibly zero) virtual projective character which is a $\mathbb{Z}$-linear combination of irreducible characters in $b$.

\medskip
Generalizing some of the results of the previous section, we may note that:

\medskip
\noindent {\bf Theorem 4.1:} \emph{ If $b$ is a block of $RG$, then $\Psi^{(b)}_{1,p,G}$  contains the character of at least one projective indecomposable $b$-module with non-negative multiplicity. In particular, if $\ell(b) = 1,$ then $\Psi^{(b)}_{1,p,G}$ is either zero or else is a character of $G$ which may be afforded by a projective $RG$-module.}

\medskip
\noindent {\bf Proof:} Let $\theta_{i}$ be the character afforded by some projective indecomposable $b$-module. Then $$\langle \theta_{i}, \Psi_{1,p,G} \rangle \geq 0, $$ since $\theta_{i}$ vanishes on all $p$-singular elements and $\Psi_{1,p,G}$ agrees with a permutation character on $p$-regular elements.

\medskip
Now $\theta_{i}$ agrees with a non-negative integer linear combination of Brauer irreducible characters in $b$ on $p$-regular elements. Hence there must be some $j$ with $$\langle \phi_{j}, \Psi_{1,p,G} \rangle \geq 0$$ (and with $\phi_{j} \in b$). Then the character $\theta_{j}$ is a character of a projective indecomposable $b$-module, and occurs with non-negative multiplicity when $\Psi_{1,p,G}$ is expressed as a $\mathbb{Z}$-linear combination of $\{ \theta_{r}: 1 \leq r \leq \ell \}.$

\medskip
\noindent {\bf Remark 4.2 : } This provides an alternative proof of part i) of Theorem 3.2. It also reproves the fact that irreducible characters in $p$-blocks of defect zero all occur with non-negative multiplicity in $\Psi_{1,p,G}.$

\section{Other $p^{\prime}$-sections}%5

\medskip
Let $y$ be a $p$-regular element of $G$, and consider the generalized character $\Psi_{1,p,C_{G}(y)}$ of $C_{G}(y),$ which vanishes on all $p$-singular elements of $C_{G}(y)$, and takes value (\emph{the number of $p$-elements of $C_{G}(y) \cap C_{G}(z)$}) for each $p$-regular $z$ in $C_{G}(y)$. This is the exact analogue of $\Psi_{1,p,G},$ but for the subgroup $C_{G}(y)$.

\medskip
Then $\Psi_{1,p,C_{G}(y)}$ agrees with 
$$\sum_{ x \in C_{G}(y)_{p}/C_{G}(y)} {\rm Ind}_{C_{G}(xy)}^{C_{G}(y)}(1)$$ on $p$-regular elements of $C_{G}(y)$ and vanishes on all $p$-singular elements of $C_{G}(y)$. Notice then that 
${\rm Ind}_{C_{G}(y)}^{G}(\Psi_{1,p,C_{G}(y)})$ agrees on $p$-regular elements of $G$ with the character afforded by the conjugation action of $G$ on the $p^{\prime}$-section of $y,$ and vanishes on all $p$-singular elements of $G$.

\medskip
Hence we see that $$\sum_{ y \in G_{p^{\prime}}/G}{\rm Ind}_{C_{G}(y)}^{G}(\Psi_{1,p,C_{G}(y)})$$ agrees on $p$-regular elements of $G$ with the character afforded by the conjugation action of $G$ on itself, and vanishes on all $p$-singular elements of $G$. Thus, by block orthogonality relations,  we have:

\medskip
\noindent {\bf Remark 5.1:} $$ \sum_{y \in G_{p^{\prime}}/G} {\rm Ind}_{C_{G}(y)}^{G}(\Psi_{1,p,C_{G}(y)})
= \sum_{i =1}^{\ell} \overline{\phi_{i}} \theta_{i}.$$  We call this character the \emph{truncated conjugation character of G}, and denote it by $\Lambda_{c,p,G}$. We denote the projective $RG$-module affording this character by $T_{c,p,G}$ and call it the \emph{truncated conjugation module}.

\medskip
Note that $\Lambda_{c,p,G}$ takes value $|C_{G}(y)|$ on each $p$-regular element $y \in G$ (and vanishes on all $p$-singular elements). Note also that $T_{c,p,G}$ is  the lift to an $RG$-module of the 
projective $\mathbb{F}G$-module $ \bigoplus_{S} {\rm Hom}(S,P(S)),$  where $S$ runs over a full set of isomorphism types of absolutely irreducible $\mathbb{F}G$-modules, and $P(S)$ denotes the projective cover of $S$ as $\mathbb{F}G$-module. 

\medskip
We remark that the multiplicity of the projective cover of the trivial module as a summand of the $\mathbb{F}G$-module $ \bigoplus_{S} {\rm Hom}(S,P(S)),$ is  $\ell(G),$ the number of simple $\mathbb{F}G$-modules, which is also the number of $p$-regular conjugacy classes of $G$.

\medskip
We note that if $Z$ is a central $p^{\prime}$-subgroup of $G$, then for each $p$-regular $y \in G$, the number of $p$-elements in $C_{G}(y)$ is the same as the number of $p$-elements of $C_{G/Z}(yZ).$
For if $x$ is a $p$-element of $G$ with $[y,x] \in yZ,$ then we have $x \in C_{G}(y)$ by a standard argument, and $x$ is the unique $p$-element in $xZ.$ Hence $\Psi_{1,p,G}$ is the inflation to $G$ of $\Psi_{1,p,G/Z}$, and $\Psi_{1,p,G}$ is afforded by a projective module if and only if 
$\Psi_{1,p,G/Z}$ is afforded by a projective module.

\medskip
Hence to prove that $\Psi_{1,p,G}$ is afforded by a projective $RG$-module, we may assume by induction that $Z(G)$ is a $p$-group. But we have already dealt with the case $O_{p}(G) \neq 1,$  so to prove that $\Psi_{1,p,G}$ is a character afforded by a projective $RG$-module, we may assume by induction that $O_{p}(G) = Z(G) = 1.$

\medskip
\noindent {\bf Example 5.2:} If we consider the special case that $G = {\rm GL}(n,p^{a})$ for some positive integer $a$, then for each $p$-regular element $y \in G$, it is clear (and well-known) that $C_{G}(y)$ is a direct product of general linear groups over fields of characteristic $p$ (allowing the possibility that a cyclic group of order $p^{r}-1$ is considered as ${\rm GL}(1,p^{r})$).  The Steinberg character ${\rm St}_{1}$ of $G$ is known to have degree 
$|G|_{p}$ and, more generally, to take value $\pm |C_{G}(y)|_{p}$ for each $p$-regular $y \in G$.

\medskip
Hence, for each $p$-regular $y \in G,$ we see that  $C_{G}(y)$ has a Steinberg character ${\rm St}_{y}$, and we have ${\rm St}_{1}(y) = \pm {\rm St}_{y}(1)$ for each $p$-regular $y \in G$. Since ${\rm St}_{1}$ is irreducible and vanishes on $p$-singular elements, we have $$|G| = \sum_{y \in G_{p^{\prime}/G} } [G:C_{G}(y)] {\rm St}_{y}(1)^{2}.$$

\medskip
Since we have \\
$|G| = \sum_{y \in G_{p^{\prime}/G} } [G:C_{G}(y)] \times $  \emph{(the number of $p$-elements of $C_{G}(y))$}, an inductive argument recovers the well-known fact  that for each $p$-regular $y \in G$, the number of $p$-elements in $C_{G}(y)$ is equal to $|C_{G}(y)|_{p}^{2}$.

\medskip
Notice then that now we have $\Psi_{1,p,G}  = {\rm St}_{1}^{2},$ so proving that $\Psi_{1,p,G}$ is afforded by a projective $RG$-module, since the Steinberg character of $G$ is afforded by a projective $RG$-module for this choice of $G$. 

\medskip
In fact, the same statement holds whenever $G$ is a finite group which has  a (split,restricted) $BN$-pair with characteristic $p$ (the number of $p$-elements of $C_{G}(y)$ is equal to $|C_{G}(y)|_{p}^{2}$ whenever $y \in G$ is $p$-regular, and consequently $\Psi_{1,p,G}$ is afforded by a projective $RG$-module, for all such $G$). We will provide a new proof of this result later.

\medskip
Returning to the case $G = {\rm GL}(n,p^{a})$, we note that we have 
$$\Lambda_{c,p,G} = \sum_{y \in G_{p^{\prime}}/G} {\rm Ind}_{C_{G}(y)}^{G} ( {\rm St}_{y}^{2}).$$

\medskip
\section{Some consequences of non-negativity}%6

\medskip
When $\phi$ is a Brauer character of $G$ (with respect to the prime $p$), we define the generalized character $\phi^{\ast}$ of $G$ via $\phi^{\ast}(xy) = \phi(y)$ whenever $x$ is a $p$-element and $y$ is a $p$-regular element with $xy =yx$. That this is a generalized character of $G$ is an easy consequence of Brauer's characterization of characters, as noted by J.A. Green in [8].

\medskip
\noindent {\bf Theorem 6.1:} \emph{Let $G$ be a finite group. Then $\Psi_{1,p,G}$ is a non-negative integer combination of characters of projective indecomposable $RG$-modules if and only if we have $\langle \phi_{i}^{\ast}, \phi_{j}^{\ast} \rangle \geq 0$ for $1 \leq i,j \leq \ell$.}

\medskip
\noindent {\bf Proof:} Let $\phi_{i}$ be an irreducible Brauer character of $G$.

\medskip
Then the coefficient of $\theta_{i}$ when $\Psi_{1,p,G}$ is expressed as a  $\mathbb{Z}$-linear combination of characters of projective indecomposable $RG$-modules is  $$\langle \Psi_{1,p,G}, \phi_{i} \rangle = \langle \Psi_{1,p,G}, \phi_{i}^{\ast} \rangle$$  (recall that $\Psi_{1,p,G}$ vanishes on $p$-singular elements, so on all elements of the $p^{\prime}$-section of $y$ outside the conjugacy class of $y$).

\medskip
Since \emph{$\Psi_{1,p,G}(y) = $ (the number of $p$-elements of $C_{G}(y)$)} for each $p$-regular $y \in G,$ we see from the definitions of $\Psi_{1,p,G}$ and $\phi_{i}^{\ast}$ that this last inner product is $\langle 1,\phi_{i}^{\ast} \rangle.$ 

\medskip
Hence $\Psi_{1,p,G}$ is the character afforded by a projective $RG$-module if and only if $$\langle \Psi_{1,p,G}, \phi_{i} \rangle  = \langle 1,  \phi_{i}^{\ast}\rangle$$ is non-negative for each absolutely irreducible Brauer character $\phi_{i}$ of $G$. Suppose then that this positivity condition holds.

\medskip
It is clear that if $\phi_{i}$ and $\phi_{j}$ are irreducible Brauer characters, then $\overline{\phi_{i}}^{\ast}\phi_{j}^{\ast}$ is a non-negative integer combination of extended (that is to say, starred) Brauer irreducible characters, so this implies that $$\langle \phi_{i}^{\ast}, \phi_{j}^{\ast} \rangle = \langle \overline{\phi_{i}}^{\ast} \phi_{j}^{\ast}, 1  \rangle \geq 0.$$

\medskip
Hence if $\Psi_{1,p,G}$ is the character of a projective $RG$-module, then $$\langle \phi_{i}^{\ast},\phi_{j}^{\ast} \rangle \geq 0$$ for all $i,j$. 

\medskip
On the other hand, if such inner products are always non-negative, then the case that $\phi_{j}$ is the trivial character already implies that $$\langle \Psi_{1,p,G}, \phi_{i} \rangle = \langle 1, \phi_{i}^{\ast} \rangle \geq 0$$ for each $i$, so that $\Psi_{1,p,G}$ is a non-negative integer combination of characters of projective indecomposable $RG$-modules.

\medskip
\noindent {\bf Remark 6.2 :} Since the generalized characters $\phi_{i}^{\ast}$ are not, in general, characters, this non-negativity condition appears to be slightly surprising.

\section{ Some remarks on $\Psi_{1,p,G}$ for groups in which the centralizer of each non-identity $p$-element has a normal $p$-complement}%7

\medskip
We suppose (throughout this section) that the centralizer of each non-identity $p$-element of $G$ has a normal $p$-complement.

\medskip
We may decompose $\Psi_{1,p,G}$ according to blocks of $RG$. We let $B_{0}$ denote the principal block of $RG$,  and we set  
$$\Psi^{(0)}_{1,p,G} = 
\sum_{ \chi \in {\rm Irr}(B_{0})} \langle \Psi_{1,p,G},  \chi \rangle  \chi .$$  Then we define $\Psi^{(1)}_{1,p,G}$ as $\Psi_{1,p,G}- \Psi^{(0)}_{1,p,G}.$

\medskip
\noindent {\bf Theorem 7.1: } \emph{$ \Psi^{(1)}_{1,p,G}$ is a character afforded by a projective $RG$-module.}

\medskip
\noindent {\bf Proof: } It suffices to prove that whenever $\phi_{i}$ is an absolutely irreducible Brauer character in a non-principal block of $RG$, then 
$$\langle \Psi_{1,p,G}, \phi_{i}\rangle \geq 0.$$

\medskip
We have  
$$\langle \Psi_{1,p,G}, \phi_{i} \rangle  = \sum_{x \in G_{p}/G} \langle {\rm Res}^{G}_{C_{G}(x)}(\phi_{i}),1 \rangle^{\prime}$$ and (since $C_{G}(x)$ has a normal $p$-complement whenever $1 \neq x$ is a $p$-element) the contribution from non-identity $p$-elements $x$ is non-negative.

\medskip
But by block orthogonality relations, we know that $\sum_{g \in G_{p^{\prime}}}\phi_{i}(g) = 0$ whenever $\phi_{i}$ is a Brauer irreducible character not in the principal block of $RG$, so the required non-negativity follows for all Brauer characters in non-principal blocks. 

\medskip 
As a consequence of 7.1, we have:

\medskip
\noindent {\bf Corollary 7.2 : } \emph{ Let $G$ be a finite group whose Sylow $p$-subgroups  are cyclic, dihedral or Klein four, or that $G$ has an Abelian Sylow $p$-subgroup $P$ such that $N_{G}(P)/O_{p^{\prime}}(C_{G}(P))$ is a Frobenius group with Frobenius kernel isomorphic to $P$.}

\medskip
\emph{Then we may write 
$\Psi_{1,p,G}  = \Psi^{(0)}_{1,p,G} + \Psi^{(1)}_{1,p,G},$ where $\Psi^{(0)}_{1,p,G}$ is a virtual projective character which is a $\mathbb{Z}$-combination of irreducible characters in the principal $p$-block, and $\Psi^{(1)}_{1,p,G}$ is the character of a projective $RG$-module which has no indecomposable summand in the principal $p$-block of $G$.} 

\medskip
\noindent {\bf Proof:} In all of the cases listed, $C_{G}(x)$ has a normal $p$-complement for each non-identity $p$-element $x \in G$.

\medskip
\noindent {\bf Theorem 7.3} \emph{ Suppose that $C_{G}(x)$ has a normal $p$-complement for each non-identity $p$-element $x \in G.$ Then $\Psi_{1,p,G}$ is a character of $G$.}

\medskip
\noindent {\bf Proof:} By 7.1, it suffices to prove that $\langle \Psi_{1,p,G}, \chi \rangle \geq 0$ for each irreducible character $\chi$ in the principal $p$-block of $G$. It is clearly sufficient to consider the case that $\chi$ is non-trivial, since (for example), we always have $\langle \Psi_{1,p,G}, 1 \rangle = 1.$

\medskip
We first note that if $\chi$ is a non-trivial irreducible character in the principal $p$-block of $G$, then we have $\sum_{g \in G} \chi(g) = 0,$ and from this (and from Brauer's Second and Third Main Theorems, since they imply that $\chi$ is constant on non-identity $p$-sections), we 
easily deduce that 
$$\frac{1}{|G|}\left( \sum_{g \in G_{p^{\prime}}} \chi(g) \right) = -\sum_{1 \neq x \in G_{p}/G} \frac{\chi(x)}{|C_{G}(x)|_{p}}.$$

\medskip
We may note as in the previous proofs that for each non-identity $p$-element $x$ of $G$, we have 
$$\frac{1}{|C_{G}(x)|} \left( \sum_{ y \in C_{G}(x)_{p^{\prime}}}\chi(y) \right)
= \frac{1}{|C_{G}(x)|_{p}} \langle {\rm Res}^{G}_{O_{p^{\prime}}(C_{G}(x))}(\chi),1 \rangle.$$

\medskip
Now let $x$ be a non-identity $p$-element of $G$. By Brauer's Second Main Theorem, $\chi(x)$ can be evaluated by considering only the irreducible constituents of ${\rm Res}^{G}_{C_{G}(x)}(\chi)$ which lie in the principal $p$-block, and hence, in particular, have $O_{p^{\prime}}(C_{G}(x))$ in their kernel. This certainly yields $|\chi(x)| \leq  \langle {\rm Res}^{G}_{O_{p^{\prime}}(C_{G}(x))}(\chi),1 \rangle,$ since $\chi(x)$ is a sum of  $\langle {\rm Res}^{G}_{O_{p^{\prime}}(C_{G}(x))}(\chi),1 \rangle$ roots of unity. More precisely, we may conclude that (for each $p$-element $x$), $$\langle {\rm Res}^{G}_{O_{p^{\prime}}(C_{G}(x))}(\chi),1 \rangle - \chi(x)$$ is a complex number with non-negative real part.

\medskip
Now we have $$\langle \Psi_{1,p,G}, \chi \rangle = \sum_{1 \neq x \in G_{p}/G}\frac{1}{|C_{G}(x)|_{p}}
\left( \langle {\rm Res}^{G}_{O_{p^{\prime}}(C_{G}(x))}(\chi),1 \rangle - \chi(x) \right).$$ This is a rational integer with non-negative real part, so is a non-negative integer. Hence we do have $\langle \Psi_{1,p,G}, \chi \rangle \geq 0,$ and $\Psi_{1,p,G}$ is a character, as $\chi$ was an arbitrary non-trivial irreducible character in the principal $p$-block.

\medskip
Notice that in fact $\langle \Psi_{1,p,G}, \chi \rangle$ is strictly positive for $\chi$ as above unless $$\chi(g) =  \langle {\rm Res}^{G}_{O_{p^{\prime}}(C_{G}(g_{p}))}(\chi),1 \rangle $$ for all $p$-singular $g \in G.$ In particular, $\langle \Psi_{1,p,G}, \chi \rangle$ is strictly positive unless $\chi$ takes non-negative rational integer values at all non-identity $p$-elements of $G$. Since $\chi$ does not vanish on all non-identity $p$-elements of $G$, we may also note that if $\langle \Psi_{1,p,G}, \chi \rangle =0,$ then there is at least one non-identity $p$-element $x$ such that $\chi(x)$ takes a constant strictly positive rational integer value on the $p$-section of $x$, and, in particular, we must have $\langle {\rm Res}^{G}_{O_{p^{\prime}}(C_{G}(x))}(\chi), 1 \rangle > 0$ for this $x$.

\medskip
\noindent {\bf Corollary 7.4:} \emph{ Suppose that $G \cong {\rm PSL}(2,q)$. Then the generalized character $\Psi_{1,p,G}$
is a character of $G$ for each prime divisor $p$ of $|G|$.}

\medskip
\noindent {\bf Proof:} There are three possible structures for the Sylow $p$-subgroup $P$ of $G$, and we note that in each of these cases, the centralizer of each non-identity $p$-element of $G$ has a normal $p$-complement. If $P$ is cyclic, then $C_{G}(x)$ has a normal $p$-complement for each  
non-identity element $x \in P.$ If $p = 2$ and $P$ is dihedral (allowing a Klein 4-group), then 
$C_{G}(x)$ has a normal $p$-complement for each non-identity element $x \in P$. If $p|q$ and $P$ is non-cyclic elementary Abelian, then $C_{G}(x)$ is a $p$-group for each non-identity element $x \in P.$

\section{Variants for sets of primes}%8

\medskip
Let $\pi$ be a set of prime divisors of the order of the finite group $G$, and let $n = |G|_{\pi},$ that is to say, the largest divisor of $|G|$ which is only divisible by primes in $\pi$. We may refer to $n$ as a \emph{Hall divisor} of $|G|$.

\medskip
Let $\Psi_{1,\pi,G}$ denote the class function of $G$ which takes value $0$ on all $\pi$-singular elements of $G$
and which takes the value (\emph{ the number of } $\pi$-\emph{elements of} $C_{G}(x)$) on all $\pi$-regular elements of $G$. 

\medskip
Then $\Psi_{1,\pi,G}$ is a generalized character of $G$, and $\Psi_{1,\pi,G}(x)$ is the number of $n$-th roots of $x$ in $G$, where $n = |G|_{\pi}.$ Furthermore, $\Psi_{1,\pi,G}$ is a $\mathbb{Z}$-linear combination of characters each induced from linear characters of Brauer elementary $\pi^{\prime}$-subgroups of $G$, all proved in manner entirely analogous to the case $\pi = \{p\}.$

\medskip
It is also the case that if we let $G_{\pi^{\prime}}/G$ denote a full set of representatives for the conjugacy classes of $\pi$-regular elements of $G$, then we have 
$$\sum_{ y \in G_{\pi^{\prime}}/G} {\rm Ind}_{C_{G}(y)}^{G}(\Psi_{1,\pi,C_{G}(y)}) = \Lambda_{c,\pi,G}$$
where $\Lambda_{c,\pi,G}$ is the generalized character of $G$ which takes value $|C_{G}(x)|$ on $\pi$-regular elements $x$, and $0$ on all $\pi$-singular elements.

\medskip
There is some computational evidence that $\Psi_{1,\pi,G}$ is always a character, but this remains unproven at present. We are thus interested in determining whether or not $\Lambda_{c,\pi,G}$ is a character, for if there is a group $G$ such that $\Lambda_{c,\pi,G}$ is not a character of $G$, then there is a $\pi$-regular element $y \in G$ such that $\Psi_{1,\pi,C_{G}(y)}$ is not a character of $C_{G}(y).$

\medskip
We recall that a $\pi$-separable finite group $G$ is a finite group whose composition factors are all either $\pi$-groups or $\pi^{\prime}$-groups. Every solvable group is $\pi$-separable for every set of primes $\pi.$

\medskip
\noindent {\bf Theorem 8.1:} \emph{ If $G$ is a $\pi$-separable group, then $\Lambda_{c,\pi,G}$ is a character of $G$,
and $\Lambda_{c,\pi,G}$ is induced from a character of a Hall $\pi^{\prime}$-subgroup of $G$.}
 
\medskip
\noindent {\bf Proof:} Since $G$ is $\pi$-separable, $G$ has a Hall $\pi^{\prime}$-subgroup $H$. By the extension of Fong's theory for $p$-solvable groups in [7] to sets of primes for $\pi$-separable groups by I.M. Isaacs in [9], we know that if $G$ has $\ell$ conjugacy classes of $\pi$-regular elements of $G$, then there are $\ell$ irreducible characters $\{\phi_{i} : 1 \leq i \leq \ell \}$ of $G/O_{\pi}(G)$ (which we inflate to irreducible characters of $G$), such that for each irreducible character $\chi$ of $G$,  the restriction $\chi_{0}$ of $\chi$ to $G_{\pi^{\prime}}$ is a non-negative integer combination of $\{ \phi_{i,0}: 1 \leq i \leq \ell\}$, and there are $\ell$ distinct irreducible characters of $H$, $\{ \alpha_{i}: 1 \leq i \leq \ell \}$ such that $$\langle {\rm Res}^{G}_{H}(\phi_{i}), \alpha_{j} \rangle = \delta_{ij}.$$ In particular, this implies that $\{\phi_{i}: 1 \leq i \leq \ell \}$ is linearly independent over $\mathbb{C}.$

\medskip
For each $i$, if we let $\theta_{j}$ denote the induced character ${\rm Ind}_{H}^{G}(\alpha_{j}),$
then we have $\langle \theta_{j}, \alpha_{i} \rangle = \delta_{i,j}.$ Notice then that $\{ \theta_{j} : 1 \leq j \leq \ell \}$ is linearly independent over $\mathbb{C}$, and forms a $\mathbb{C}$-basis for the space of complex-valued class functions of $G$ which vanish on all $\pi$-singular elements of $G$.

\medskip
For if $\Lambda$ is any class function of $G$ which vanishes on all $\pi$-singular elements of $G$, then 
we have $$\Lambda  = \sum_{i=1}^{\ell} \langle \Lambda, \phi_{i} \rangle \theta_{i}.$$

\medskip
In particular, $\Lambda$ is a generalized character if and only if $\langle \Lambda, \phi_{i} \rangle  \in \mathbb{Z}$ for each $i$, and $\Lambda$ is a character if and only if $\langle \Lambda, \phi_{i} \rangle  \in \mathbb{N} \cup \{0\}$ for each $i$. Notice, in particular, that $\Lambda$ is a character if and only if $\Lambda$ is induced from a character of $H$. 

\medskip
So far, all is analogous to Fong's theory for $p$-solvable groups, and is as developed by I.M. Isaacs for $\pi$-separable groups.

\medskip
Here, the $\phi_{i}$ are analogous to the Brauer characters in the usual modular theory. Now let $y$ be a $\pi$-regular element of $G$, and let $\chi_{y}$ denote the characteristic function of the conjugacy class of 
$y$. Then $|C_{G}(y)|\chi_{y}$ is a $\mathbb{C}$-linear combination of irreducible characters of $G$, and certainly vanishes on $\pi$-singular elements of $G$, and we see easily using inner products as above that  
$$|C_{G}(y)|\chi_{y} = \sum_{i=1}^{\ell} \phi_{i}(y^{-1})\theta_{i}.$$

\medskip
Since each $\theta_{i}$ vanishes on $\pi$-singular elements, we conclude that the truncated conjugation class function $\Lambda_{c,\pi,G}$ coincides with the character $$\sum_{i=1}^{\ell} \overline{\phi}_{i}\theta_{i},$$ and is induced from a character of $H$ (in fact from a non-negative integer combination of $\{ \alpha_{i}: 1 \leq i \leq \ell \}$), and this linear combination of the $\alpha_{i}$ is unique subject to inducing to $\Lambda_{c,\pi,G}$).

\medskip
Now let us note that the coefficient of $\theta_{i}$ in $\Lambda_{c,\pi,G}$ is $\langle \Lambda_{c,\pi,G}, \phi_{i} \rangle,$ and this is clearly $\sum_{j=1}^{\ell} \phi_{i}(y_{j})$, where $\{ y_{j}: 1 \leq j \leq \ell \}$ is a set of representatives for the conjugacy classes of $\pi$-regular elements of $G$ (for convenience, we choose these to be closed under inversion). Note, in particular, that this sum is always a non-negative integer (cf. a theorem of L.Solomon
(in [16]) about complex irreducible characters).

\medskip
\noindent {\bf Remark 8.2:} C. Schroeder has informed us that he has independently obtained in unpublished work a different proof that the generalized character we call $\Lambda_{c,\pi, G}$ is a character when $G$ is $\pi$-separable.

\medskip
This seems an opportune moment to note that R. Boltje's theory of Explicit Brauer Induction (see [1]) (which is a form of Brauer's Induction Theorem commuting with restriction of characters) may be relevant to the question of whether $\Psi_{1,\pi,G}$ is a character of $G$ in general.

\medskip
Boltje's explicit Brauer induction theorem is compatible with Adams operations in a remarkable fashion. It guarantees that for each irreducible character 
$\chi$ of $G$, we may write $$\chi = \sum_{(H,\lambda)/G} a_{H,\lambda, \chi}{\rm Ind}_{H}^{G}(\lambda)$$ with each $a_{H,\lambda,\chi} \in \mathbb{Z},$ each $\lambda$  a linear character of the subgroup $H$ of $G$, and furthermore, we also have (for every integer $n$), $$\chi^{(n)} = \sum_{(H,\lambda)/G} a_{H,\lambda, \chi}{\rm Ind}_{H}^{G}(\lambda^{n}),$$ where $\chi^{(n)}(g) = \chi(g^{n})$ for all $g \in G$.

\medskip
Let us reconsider $\Psi_{1,p,G}$ in this light. In Section 1, we saw that $\Psi_{1,p,G}$ is a character of $G$ if and only if $\langle \chi^{(q)}, 1 \rangle \geq 0$ for each irreducible character $\chi$ of $G$, where $q$ is a power of $p$ with $q \geq |G|_{p}$ and $q \equiv 1$ (mod $|G|_{p^{\prime}}$).

\medskip
If we write, as above, $$\chi = \sum_{(H,\lambda)/G} a_{H,\lambda, \chi}{\rm Ind}_{H}^{G}(\lambda)$$ such that we also have $$\chi^{(q)} = \sum_{(H,\lambda)/G} a_{H,\lambda, \chi}{\rm Ind}_{H}^{G}(\lambda^{q}),$$ then we see easily that 
$$\langle \chi^{(q)}, 1 \rangle = \sum_{(H,\lambda): \lambda^{q} = 1 } a_{H,\lambda, \chi}.$$ Using similar arguments for $\pi$ we may deduce:

\medskip
\noindent {\bf Theorem 8.3:} \emph{i) Let $G$ be a finite group. Let $p$ be a prime divisor of $|G|$ and let $S$ be a Sylow $p$-subgroup of $G$.
Write (for each irreducible character $\chi$ of $G$, as we may) $$\chi = \sum_{(H,\lambda)/G} a_{H,\lambda, \chi}{\rm Ind}_{H}^{G}(\lambda),$$ where each $a_{H,\lambda,\chi} \in \mathbb{Z}$ and $(H,\lambda)$ ranges over $G$-conjugacy classes of pairs such that $H$ is a subgroup of $G$ and $\lambda$ is a linear character of $H$, where we also have (for every integer $n$), 
$$\chi^{(n)} = \sum_{(H,\lambda)/G} a_{H,\lambda,\chi}{\rm Ind}_{H}^{G}(\lambda^{n}),$$ where $\chi^{(n)}(g) = \chi(g^{n})$ for each $g \in G.$ Then $\Psi_{1,p,G}$ is a character of $G$ if and only if $$\sum_{(H,\lambda,\chi): \lambda^{|S|} = 1 } a_{H,\lambda,\chi} \geq 0 $$ for each irreducible character $\chi$ of $G$.}

\medskip
\emph{ii) Let $G$ be a finite group. Let $\pi$ be a set of primes and let $m$ be the $\pi$-part of $|G|$. Write (for each irreducible character $\chi$ of $G$, as we may) $$\chi = \sum_{(H,\lambda)/G} a_{H,\lambda, \chi}{\rm Ind}_{H}^{G}(\lambda),$$ where each $a_{H,\lambda,\chi} \in \mathbb{Z}$ and $(H,\lambda)$ ranges over $G$-conjugacy classes of pairs such that $H$ is a subgroup of $G$ and $\lambda$ is a linear character of $H$, where we also have 
$$\chi^{(n)} = \sum_{(H,\lambda)/G} a_{H,\lambda,\chi}{\rm Ind}_{H}^{G}(\lambda^{n}),$$ where $\chi^{(n)}(g) = \chi(g^{n})$ for each $g \in G.$ Then $\Psi_{1,\pi,G}$ is a character of $G$ if and only if $$\sum_{(H,\lambda): \lambda^{m} = 1 } a_{H,\lambda,\chi} \geq 0 $$ for each irreducible character $\chi$
of $G$.}

\newpage
\section{The prime to the first power case, and groups with a strongly $p$-embedded subgroup}%9

\medskip
Let $G$ be a finite group with Sylow $p$-subgroup $P$ of order $p$ for $p$ a fixed prime. Let $H = N_{G}(P)$ and let $T$ be a Hall $p^{\prime}$-subgroup of $H$. Then from Theorem 3.2, we know that $\Psi_{1,p,H}$ is the character afforded by the projective cover of the trivial module, which is ${\rm Ind}_{T}^{H}(1)$. 

\medskip
We also know that $$\Psi_{1,p,H} = \chi_{r,H} + \sum_{ 1 \neq x \in P /H} {\rm Ind}_{C_{H}(x)}^{H}(\chi_{r,C_{H}(x)}),$$ where $\chi_{r,H}$ is the characteristic function of the set of $p$-regular elements of $H,$ etc.

\medskip
Note that $H$ controls the conjugacy of elements of $P$ in $G$, and notice also that
 $${\rm Ind}_{H}^{G}(\chi_{r,H}) - \chi_{r,G}
= {\rm Ind}_{H}^{G}(1) - 1.$$

\medskip
Now we have $$\Psi_{1,p,G} = \chi_{r,G} + \sum_{ 1 \neq x \in P/G} {\rm Ind}_{C_{G}(x)}^{G}( \chi_{r,C_{G}(x)}).$$

\medskip
Now we see that $${\rm Ind}_{H}^{G}(\Psi_{1,p,H}) - \Psi_{1,p,G} = {\rm Ind}_{H}^{G}(1)-1 . $$

\medskip
Hence we have $$\Psi_{1,p,G} = {\rm Ind}_{H}^{G}(\Psi_{1,p,H}) - {\rm Ind}_{H}^{G}(1)  + 1.$$ 

\medskip
Since we know that $$\Psi_{1,p,H} = {\rm Ind}_{T}^{H}(1),$$ we  have now proved that in this case, we have:

\medskip
\noindent {\bf Theorem 9.1:} \emph{ If the group $G$ has a Sylow $p$-subgroup $P$ of order $p$, then we have 
$$\Psi_{1,p,G} = 1 + {\rm Ind}_{T}^{G}(1) - {\rm Ind}_{H}^{G}(1),$$ where $H = N_{G}(P)$ and $T$ is a Hall $p^{\prime}$-subgroup of $H$.}

\medskip
\noindent {\bf Remark 9.2:} We note that in Theorem 9.1, the irreducible constituents of $\Psi_{1,p,G}$ are among the irreducible constituents of ${\rm Ind}_{T}^{G}(1).$

\medskip
We claim next that $\Psi_{1,p,G}$ is the character afforded by some projective $RG$ module. To see this, it suffices to prove that $\Psi_{1,p,G}$ agrees on $p$-regular elements with the  Brauer character of a projective $\mathbb{F}G$-module.

\medskip
Let $P_{1}$ denote the projective cover of the trivial $\mathbb{F}H$-module $\mathbb{F},$ whose Brauer character agrees with $\Psi_{1,p,H}$ on $p$-regular elements. Let us write $${\rm Ind}_{H}^{G}(\mathbb{F})  =  \mathbb{F} \oplus Q,$$ where $Q$ is projective. We need to prove that ${\rm Ind}_{H}^{G}(P_{1})$ has an indecomposable summand isomorphic to $Q$. Since projective $\mathbb{F}G$-modules are injective, it suffices to prove that ${\rm Ind}_{H}^{G}(P_{1})$ has a submodule isomorphic to $Q.$ But ${\rm Ind}_{H}^{G}(\mathbb{F})$ is certainly isomorphic to a submodule of ${\rm Ind}_{H}^{G}(P_{1}),$ so this is clear.
  
\medskip
\noindent {\bf Theorem 9.3:} \emph{Suppose that $G$ is a finite group which has a (proper) strongly $p$-embedded subgroup $H$. Suppose further that $\Psi_{1,p,H}$ may be afforded by a projective $RH$-module $M$. Then $\Psi_{1,p,G}$ may be afforded by a projective $RG$-module.}

\medskip
\noindent {\bf Proof:} We mimic the proof in the case that the Sylow $p$-subgroup of $G$ had order $p$. Just as in the preamble to 9.1, we see 
$${\rm Ind}_{H}^{G}(\Psi_{1,p,H}) - \Psi_{1,p,G} = {\rm Ind}_{H}^{G}(1)-1 , $$ so that $$\Psi_{1,p,G} = 1 + {\rm Ind}_{H}^{G}(\Psi_{1,p,H}) -
{\rm Ind}_{H}^{G}(1).$$ Since $\Psi_{1,p,H}$ is assumed to be a character, and contains the trivial character with multiplicity one, it is clear that 
$\Psi_{1,p,G}$ is a character (still containing the trivial character with multiplicity one). 

\medskip 
Now let $M$ be a projective $\mathbb{F}H$-module  whose Brauer character agrees with $\Psi_{1,p,H}$ on $p$-regular elements of $H$. Since $H$ is strongly $p$-embedded in $G$,  we may write $${\rm Ind}_{H}^{G}(\mathbb{F})  =  \mathbb{F} \oplus Q,$$ where $Q$ is a projective $\mathbb{F}G$-module.  

\medskip
It suffices to prove that, in the Green ring of $\mathbb{F}G$-modules, $$\mathbb{F}  + {\rm Ind}_{H}^{G}(M) - {\rm Ind}_{H}^{G}(\mathbb{F})$$ represents a genuine projective $\mathbb{F}G$-module. For this, it is enough to show that ${\rm Ind}_{H}^{G}(M)$ has a direct summand isomorphic to $Q$. Now $M$ has a direct summand which is the projective cover of the trivial $\mathbb{F}H$-module. Hence the socle of $M$ has a trivial summand and ${\rm Ind}_{H}^{G}(\mathbb{F})$ is certainly isomorphic to a submodule of  ${\rm Ind}_{H}^{G}(M)$. In particular, $Q$ is isomorphic to a direct summand of ${\rm Ind}_{H}^{G}(M)$ (since $Q$ is projective, and hence also injective), as required.

\medskip
\noindent {\bf Remark 9.4:} In fact, $M$ has at least two trivial composition factors since the projective cover of the trivial $\mathbb{F}H$-module has trivial head and trivial socle, so that $Q$ occurs with multiplicity at least two as a summand of ${\rm Ind}_{H}^{G}(M)$, and the lift of $Q$ to a projective $RG$-module occurs as a direct summand of the projective $RG$-module affording $\Psi_{1,p,G}.$ Furthermore, it is easy to check that the projective cover of the trivial $\mathbb{F}G$-module  is not isomorphic to a direct summand of $Q$.

\medskip
\noindent {\bf Theorem 9.5:} \emph{ If $G$ has a cyclic Sylow $p$-subgroup $S$, then $\Psi_{1,p,G}$ is a character afforded by a projective $RG$-module.}

\medskip
\noindent {\bf Proof:} We proceed by induction on $|S|$, the result being clear when $|S| = 1, $ so suppose that $|S| >1$  and the result is true for groups $G$ with a Sylow $p$-subgroup of order less than  $|S|$. If $O_{p}(G) \neq 1,$ then by induction $\Psi_{1,p,G/O_{p}(G)}$ is afforded by a projective $RG/O_{p}(G)$-module $X$, and then $\Psi_{1,p,G}$ is afforded by the projective cover of (the inflation of) $X$ as $RG$-module. Hence we may suppose that $O_{p}(G) = 1.$ 

\medskip
But then $N = N_{G}(\Omega_{1}(P))$ is a proper strongly $p$-embedded subgroup of $G$, where $\Omega_{1}(P)$ is the unique subgroup of $P$ of order $p$, and furthermore, $\Psi_{1,p,N}$ is afforded by a projective $RN$-module by the argument of the previous paragraph. By Theorem 9.3, $\Psi_{1,p,G}$ is afforded by a projective $RG$-module.

\medskip
We also mention another case of interest, which is a special case of Theorem 9.3, after noting that $$\Psi_{1,p,H} = {\rm Ind}_{T}^{H}(1)$$ by part iii) of Theorem 3.2.

\medskip
\noindent {\bf Corollary 9.6:} \emph{Suppose that $G$ has a $TI$-Sylow $p$-subgroup $P$. Let $H = N_{G}(P)$ and $T$ be a Hall $p^{\prime}$-subgroup of $H$. Then $$\Psi_{1,p,G}=1 +  {\rm Ind}_{T}^{G}(1) - {\rm Ind}_{H}^{G}(1)$$ is a character which may be afforded by a projective $RG$-module.}

\medskip
\section{More on the Truncated Conjugation Module}%10

\medskip
We recall that the $T_{c,p,G},$ the truncated conjugation module, is the projective $RG$-module whose character takes value $0$ on each $p$-singular element of $G$, and value $|C_{G}(x)|$ on each $p$-regular element of $G$. We propose the following:

\medskip
\noindent {\bf Conjecture A:} \emph{ For each $p^{\prime}$-element $x$ of $G$, the $R$-valued class function $\Psi_{x,p,G}$ which agrees on $p$-regular elements with the character afforded by  the conjugation action on the $p^{\prime}$-section of $x$ in $G$, and which vanishes on all $p$-singular elements, is  a character, and may be afforded by a projective $RG$-module $P_{x,p,G}$. Consequently, 
$$T_{c,p,G}= \bigoplus_{ x \in G_{p^{\prime}/G}} P_{x,p,G}.$$}

\medskip
\noindent {\bf Lemma 10.1:} \emph{ To prove Conjecture A, it suffices to treat (for all $G$) the case that $x = 1_{G}.$}

\medskip
\noindent {\bf Proof:} Suppose that for every finite group $H$, we know that the class function $\Psi_{1,p,H}$ which agrees with the permutation character
$$\sum_{y \in H_{p}/H}  {\rm Ind}_{C_{H}(y)}^{H}(1)$$ on $p$-regular elements, and vanishes elsewhere, is a character, and is afforded by a projective $RH$-module.  We claim that Conjecture A holds for $G$.

\medskip
Suppose first that $G$ contains a non-trivial central $p$-regular element $z.$ Then, as noted at the end of Remark 5.1, $$\Psi_{z,p,G}  = \Psi_{1,p,G}$$ and $\Psi_{1,p,G}$ is the inflation to $G$ of $$\Psi_{1,p,G/\langle z \rangle}.$$ By hypothesis, $\Psi_{1,p,G/\langle z \rangle}$ is afforded by a projective $RG/\langle z \rangle$-module, which may be regarded as (ie inflates to) a projective $RG$-module.

\medskip
Now let $x$ be any non-central $p$-regular element of $G$. Then Conjecture A holds for $C_{G}(x).$ The class function $\Psi_{x,p,G}$ is clearly induced from $\Psi_{x,p,C_{G}(x)}.$ This last class function is afforded by the projective $RC_{G}(x)$-module $P_{x,p,C_{G}(x)},$ arguing as in the paragraph above. Hence $\Psi_{x,p,G}$ is a character and is afforded by the projective $RG$-module ${\rm Ind}_{C_{G}(x)}^{G}(P_{x,p,C_{G}(x)}).$ Thus Conjecture A holds for $G.$

\medskip
Hence proving Conjecture A has been reduced to proving:

\medskip
\noindent {\bf Conjecture B:} \emph{ The class function $\Psi_{1,p,G}$ which takes value ({\rm the number of} $p$-{\rm elements of} $C_{G}(x)$) on each $p$-regular $x \in G,$  and which vanishes on all $p$-singular elements of $G,$ is a character of $G$ which may be afforded by a projective $RG$-module.}

\medskip
\noindent {\bf Corollary 10.2:} \emph{ Conjecture B holds for $G = {\rm PSL}(2,q)$ for every odd prime $p$, and also holds for $p=2$ if $q$ is even.}

\medskip
\noindent {\bf Proof:} Let $H = N_{G}(P),$ where $P$ is a Sylow $p$-subgroup of $G$. If $p$ is odd, then either $P$ is cyclic, or else $q$ is a power of $p$ and $P$ is elementary Abelian. If $q$ is even and $p = 2,$ then $P$ is elementary Abelian, and $H$ is strongly (2)-embedded. By Theorem 3.2, $\Psi_{1,p,H}$ is a character afforded by a projective $RH$-module  in all cases. Hence the result follows by Theorem 9.3 and Theorem 9.5.

\medskip
\noindent {\bf Remark 10.3:}  It is true that $\Psi_{1,2,G}$ is a character whenever $G = {\rm PSL}(2,q)$,  and that the contribution to $\Psi_{1,2,G}$ from non-principal $2$-blocks is afforded by a projective module, as noted in 7.1 and 7.3. 

\medskip
We note here the following alternative direct proof that $\Psi_{1,p,G}$ is a generalized character.

\medskip
\noindent {\bf Lemma 10.4:} \emph{The class function $\Psi_{1,p,G}$ is a generalized character, and is afforded by a $\mathbb{Z}$-linear combination of projective $RG$-modules (in the Green ring for $RG$).}

\medskip
\noindent {\bf Proof:} We know that $\Lambda_{c,p,G}$ is the character afforded by the projective $RG$-module $T_{c,p,G}.$ We also know that $$\Lambda_{c,p,G} = \sum_{y \in G_{p^{\prime}}/G} {\rm Ind}_{C_{G}(y)}^{G}( \Psi_{1,p,C_{G}(y)/\langle y \rangle}),$$ so the result follows by induction
(by assuming that the result is true for groups of order less than $|G|$, we may suppose that  $\Psi_{1,p,C_{G}(y)/\langle y \rangle}$ is a generalized character whenever $1 \neq y \in G_{p^{\prime}}$ (afforded by a virtual projective). We also know that $\Lambda_{c,p,G}$ is a character which may be afforded by a projective $RG$-module).

\section{Local control of the permutation action on non-identity $p$-elements}%11

\medskip
Let $\mathcal{S}_{p}(G)$ denote the simplicial complex associated to the poset of non-trivial $p$-subgroups of $G$. As usual, when $\sigma$ is an element of $\mathcal{S}_{p}(G)$, we denote the number of non-trivial $p$-subgroups in $\sigma$ by $|\sigma|.$

\medskip
The first result of this section illustrates that the conjugation action of $G$ on its non-identity $p$-elements is $p$-locally controlled in a precise sense. This is relevant to the study of $\Psi_{1,p,G},$ because $\Psi_{1,p,G} - 1 $ agrees on $p$-regular elements with the character afforded by the conjugation action of $G$ on its non-identity $p$-elements.

\medskip
\noindent {\bf Theorem 11.1:} \emph{ Let $G$ be a finite group, and $R$ be as before. For $H$ a subgroup of $G$, let $X_{p,H}$ denote the permutation module of $H$ afforded by the conjugation action on its non-identity $p$-elements. Then in the Green ring of $RG$, we have $$\sum_{ \sigma \in \mathcal{S}_{p}(G)/G} (-1)^{|\sigma|}{\rm Ind}_{G_{\sigma}}^{G}(X_{p,G_{\sigma}}) = 0.$$ In particular, the virtual character $$\sum_{ \sigma \in \mathcal{S}_{p}(G)/G} (-1)^{|\sigma|}{\rm Ind}_{G_{\sigma}}^{G}(\Psi_{1,p,G_{\sigma}}-1 ) $$ vanishes on all $p$-regular elements of $G$.}

\medskip
\noindent {\bf Proof:} We work instead with the simplicial complex $\mathcal{S}^{+}_{p}(G)$ associated to the poset of ordered pairs $(Q,x)$  whose $1$-simplices are ordered pairs $(Q,x)$ with $Q$ a non-trivial $p$-subgroup of $G$ and $x$ a non-identity $p$-element of $G$ normalizing $Q$. We decree that $(Q,x) \leq (P,y)$ if and only if $x = y$ and $Q \leq P.$

\medskip
Then it is clear that the alternating sum of permutation modules in the statement of the Theorem is equal to 
$$\sum_{ \tau  \in \mathcal{S}^{+}_{p}(G)/G} (-1)^{|\tau|}{\rm Ind}_{G_{\tau}}^{G}(R)$$ in the Green ring for $RG$, where we include the empty chain, and consider it to have stabilizer $G$ and length zero. Notice that the empty chain of $\mathcal{S}_{p}(G)$ contributes 
$$ \sum_{1 \neq x \in G_{p}/G}{\rm Ind}_{C_{G}(x)}^{G}(R)$$ to the original alternating sum of modules, and that this is cancelled by the contribution from the length one chains in $\{ (\langle x \rangle,x) : 1 \neq x \in G_{p}/G \}.$

\medskip
For the proof, we employ a familiar cancellation argument. For a given $p$-element $x \in G^{\#},$ we may reduce the contribution from chains of elements with second component $x$ to that from those chains of elements $(Q,x)$ with $[Q,x]=1.$ If $$\tau = (Q_{1},x) < (Q_{2},x) < \ldots (Q_{n},x)$$ is a chain not of this form, choose $r$ minimal  such that $[Q_{r},x] \neq 1$ (possibly $r = 1$) and form a chain $\tau^{\ast}$ by inserting $C_{Q_{r}}(x)$ between $Q_{r-1}$ and $Q_{r}$ (or at the beginning if $r = 1$) in case $Q_{r-1} \neq C_{Q_{r}}(x)$, while if $Q_{r-1} = C_{Q_{r}}(x),$ we form $\tau^{\ast}$ by deleting $Q_{r-1}$ from $\tau.$ Then it is easy to check that $\tau^{\ast \ast} = \tau,$ and that $G_{\tau^{\ast}} = G_{\tau},$ so that the contributions from  $\tau$ and $\tau^{\ast}$ cancel.

\medskip
Hence we now only need to consider (all) chains of elements $(Q,x) \in \mathcal{S}^{+}_{p}(G)$ such that $[Q,x] = 1.$ For each non-identity $p$-element $x,$ we may delete the contribution from remaining chains with second component $x$ by a similar pairing argument to that above, using respectively 
the monotonic maps $$(Q,x)  \to (Q\langle x \rangle, x)$$ and $$(Q\langle x \rangle,x) \to (\langle x \rangle,x).$$

\medskip
To expand: if we take a representative $x$ for a class of non-identity $p$-elements of $G$, and we consider a chain $$\tau = (Q_{1},x) < (Q_{2},x) < \ldots <  (Q_{n},x)$$ where each $Q_{i}$ is a subgroup of $C_{G}(x)$ containing $x$, then we may cancel the contribution to our alternating sum from the chain $\tau$ whenever $Q_{1} > \langle x \rangle$ using the chain $$(\langle x \rangle , x) < (Q_{1},x) < \ldots < (Q_{n},x)$$ in the case that $Q_{1}$ strictly contains $x$ , while if $Q_{1} = \langle x \rangle$  and $n >1$ we cancel using the chain $$(Q_{2},x) < \ldots  < (Q_{n},x).$$  If $Q_{1} = \langle x \rangle$ and $n = 1$ then we cancel the contribution from the singleton chain $(\langle x \rangle, x)$  using the pair $(1,x)$ as discussed at the beginning of the proof.

\medskip
Before our next result, we recall that $$\sum_{ \sigma \in \mathcal{S}_{p}(G)/G} (-1)^{|\sigma|} {\rm Ind}_{G_{\sigma}}^{G}(R) = {\rm St}_{p}(G),$$ is the Steinberg (virtual) module for $G$, as defined (up to sign convention) by P.J.Webb (see, e.g. [17]). This was proved by Webb to be a virtually projective $RG$-module.

\medskip
\noindent {\bf Corollary 11.2 :} \emph{ Let $P_{1,p,G}$ denote the (unique) virtual projective module in the Green ring of $RG$ affording the virtual character $\Psi_{1,p,G}$. Then we have $$\sum_{ \sigma \in \mathcal{S}_{p}(G)/G} (-1)^{|\sigma|} {\rm Ind}_{G_{\sigma}}^{G}(P_{1,p,G_{\sigma}} -R) = 0$$
in the Green ring for $RG$. In particular, $$\sum_{ \sigma \in \mathcal{S}_{p}(G)/G} (-1)^{|\sigma|} {\rm Ind}_{G_{\sigma}}^{G}(P_{1,p,G_{\sigma}}) = {\rm St}_{p}(G).$$}

\medskip
\noindent {\bf Proof:}  The virtual module in the statement of the Corollary is a virtual projective $RG$-module, as each $P_{1,p,G_{\sigma}}$ is a virtual projective. Hence to prove that it is the zero module, it is sufficient (by the non-singularity of the Cartan matrix) to prove that the alternating sum given affords the zero character. However, since it is a virtual projective, its virtual character vanishes on all $p$-singular elements. Also, as noted in Theorem 11.1, its virtual character vanishes on all $p$-regular elements, so the proof of the Corollary is complete.

\medskip
\noindent {\bf Theorem 11.3:} \emph{ In the Green ring for $RG$, the virtual projective module $P_{1,p,G}$ is equal to $$\sum_{Q/G} {\rm Ind}^{G}_{N_{G}(Q)}(P[{\rm St}_{p}\left(N_{G}(Q)/Q) \right)]),$$ where $Q/G$ is a set of representatives for the $G$-conjugacy classes of $p$-subgroups of $G$ such that $Q = O_{p}(N_{G}(Q)),$ and where $P[M]$ indicates that we are taking the (well-defined) projective cover of a virtual projective  $RN_{G}(Q)/Q$-module $M$ as (virtual) $RN_{G}(Q)$-module.}

\medskip
\noindent {\bf Proof:} We work with the complex $\mathcal{N}_{p}(G),$ where the simplices consist of chains of (strictly increasing) mutually normalizing $p$-subgroups, which  is interchangeable with $\mathcal{S}_{p}(G)$ in all calculations needed here.

\medskip
We first note that (in the Green ring for $RG$), from Corollary 11.2, we may write  
$$P_{1,p,G} - {\rm St}_{p}(G) = \sum_{1 \neq Q/G} \left(  \sum_{ \tau \in \mathcal{N}_{p}(N_{G}(Q)/Q)} (-1)^{|\tau|}{\rm Ind}_{N_{G}(Q)_{\tau}}^{G}(P_{1,p,N_{G}(Q)_{\tau}}) \right).$$ We note also that by Theorem 3.2 iii), $P_{1,p,N_{G}(Q)_{\tau}}$ is the projective cover as $RN_{G}(Q)_{\tau}$-module 
of the projective $R(N_{G}(Q)/Q)_{\tau})$-module $P_{1,p,(N_{G}(Q)/Q)_{\tau}}.$

\medskip
By Corollary 11.2, applied within $N_{G}(Q)/Q$ for each non-trivial $p$-subgroup $Q$ of $G$ (up to conjugacy), we may rewrite this last equation as 
$$P_{1,p,G} = {\rm St}_{p}(G) + 
\sum_{1 \neq Q/G} {\rm Ind}_{N_{G}(Q)}^{G} \left( P\left[{\rm St}_{p}(N_{G}(Q)/Q)\right]\right),$$ as claimed.

\medskip
\noindent {\bf Corollary 11.4:} \emph{In the Green ring of $RG$, we have $$P_{1,p,G} - R =    
\sum_{1 \neq Q/G} {\rm Ind}_{N_{G}(Q)}^{G} \left( P\left[{\rm St}_{p}(N_{G}(Q)/Q)\right] - {\rm St}_{p}(N_{G}(Q)/Q)\right).$$}

\medskip
\noindent{\bf Proof:} This follows from the fact that by Webb's inversion formula for ${\rm St}_{p}(G),$ (with suitable adjustment for sign conventions),
we have $$R = \sum_{Q/G}{\rm Ind}_{N_{G}(Q)}^{G}\left({\rm St}_{p}(N_{G}(Q)/Q)\right).$$

\medskip
\noindent {\bf Remark 11.5:} What we have shown is that the virtual projective module $$P_{1,p,G} = {\rm St}_{p}(G) + 
 \sum_{1 \neq Q/G} {\rm Ind}_{N_{G}(Q)}^{G} \left( P\left[{\rm St}_{p}(N_{G}(Q)/Q)\right]\right)$$ is the unique virtual projective in the Green ring for $RG$ which affords the virtual projective character ${\Psi}_{1,p,G}$ with\\ $\Psi_{1,p,G}(y) =$ \emph{(the number of $p$-elements of $C_{G}(y)$)} for each $p$-regular $y \in G.$

\medskip
In particular, $P_{1,p,G}$ affords a virtual character taking value $1$ on all $p$-regular elements $y$ with $p \not | |C_{G}(y)|.$ It may be instructive to verify this last fact directly. This follows from 11.4, after noting that for any non-trivial $p$-subgroup $Q$ of $G$, the Brauer characters afforded by $P\left[{\rm St}_{p}(N_{G}(Q)/Q)\right]$ and ${\rm St}_{p}(N_{G}(Q)/Q)$ agree on $p$-regular elements of defect zero of $N_{G}(Q)$, since the former Brauer character takes value $|C_{Q}(y)|$ times the value of the latter Brauer character at $y$ whenever $y \in N_{G}(Q)$ is $p$-regular (by the result of W.F. Reynolds mentioned at the end of Section 3).

\medskip
We may also give an alternative self-contained proof that $$P_{1,p,G} = {\rm St}_{p}(G) \otimes {\rm St}_{p}(G) = {\rm St} \otimes {\rm St},$$ in the case that $G$ has a (split, restricted) $BN$-pair of characteristic $p$, where ${\rm St}$ denotes the usual Steinberg module for $G$ in that case. First we need to discuss the relationship between ${\rm St}_{p}(G)$ and ${\rm St}$ in that case. It was proved by Curtis (in [4]) that the Steinberg character for $G$ is equal to $${\sum}_{ J \subseteq I} (-1)^{|J|} {\rm Ind}_{P_{J}}^{G}(1),$$ where $I$ is the set of generating reflections for the Weyl group $W$ of $G$.

\medskip
Letting $R$ denote the trivial module, this was later strengthened by P.J. Webb to the statement that in the Green ring for $RG$, we have $${\rm St} = {\sum}_{ J \subseteq I} (-1)^{|J|} {\rm Ind}_{P_{J}}^{G}(R).$$ Note that ${\rm St}$ is a genuine $RG$-module, not just a virtual module. We also recall, that $P_{\emptyset} = B,$ the Borel subgroup of $G$.  Also, Webb proved that (again in the Green ring for $RG$), we have 
$${\rm St}_{p}(G)  = \pm {\rm St}$$ when $G$ has a characteristic $p$ $BN$-pair. With the sign convention we have adopted, we see that we have 
$${\rm St}_{p}(G) = {\sum}_{ J \subseteq I} (-1)^{|I \setminus J|} {\rm Ind}_{P_{J}}^{G}(R).$$ 

\medskip
For the benefit of the reader, we clarify what we mean here by saying that $G$ has a $BN$-pair of characteristic $p$. We mean that $G$ has a split restricted $BN$-pair, and that whenever $i$ is an element of the indexing set $I$ for the generating reflections of the Weyl group $W$ for $G$, and $P_{i}$ is the associated parabolic subgroup, then $P_{i}/O_{p}(P_{i})$ has a split $BN$-pair of rank one with point stabilizer (Sylow $p$-normalizer) $B/O_{p}(P_{i})$, where $O_{p}(B)/O_{p}(P_{i})$ acts regularly on the other Sylow $p$-subgroups of $P_{i}/O_{p}(P_{i})$. 

\medskip
We now revert to standard notation for parabolic subgroups of $G$ and their unipotent radicals. At the character theoretic level, it was proved in Curtis (in [4]) that if we denote the Steinberg character of $G$ by $\chi_{I}$ and the Steinberg character of $P_{J}/U_{J}$ by $\chi_{J}$, then we have $\langle \chi_{I}, {\rm Ind}_{B}^{G}(1) \rangle = 1 $ while $$\langle \chi_{I}, {\rm Ind}_{P_{J}}^{G}(1) \rangle = 0$$ for $J \neq \emptyset.$ This readily implies ( using induction on $|I|$ and the corresponding alternating sum expression for  the character $\chi_{J}$  afforded by ${\rm St}(P_{J}/U_{J}$)),  that we have $$\langle \chi_{I}, {\rm Ind}_{P_{J}}^{G}(\chi_{J}) \rangle = 1$$ for all  $J \subseteq I.$ Using Webb's results above to translate to the Green ring, and Theorem 3 of Robinson [14], we may conclude that  $${\rm Res}^{G}_{P_{J}}({\rm St}_{G}) = P[{\rm St}_{P_{J}/U_{J}}]$$ for each $J \subseteq I,$ where the last projective cover is as an $RP_{J}$-module (strictly speaking, we use the reciprocity provided by Theorem 3 of [14] for projective simple $\mathbb{F}G$-modules and projective simple $\mathbb{F}P_{J}/U_{J}$-modules, and lift the relevant projective modules to modules over $R$).

\medskip
Now we know that ${\rm St}_{p}(N_{X}(Q)/Q)=0$ whenever $Q$ is a $p$-subgroup of a finite group $X$ with $Q \neq O_{p}(N_{X}(Q)).$ For $G$ as presently under consideration, we know that if $Q = O_{p}(N_{G}(Q))$, then $N_{G}(Q)$ is conjugate to $P_{J}$ for some $J \subseteq I$ and ${\rm St}_{p}(N_{G}(Q)/Q) \neq 0,$ so we may suppose that $Q = U_{J}$ since we deal with $p$-subgroups of $G$ up to $G$-conjugacy.

\medskip
Now we can prove:

\medskip
\noindent {\bf Theorem 11.6:} \emph{ Suppose that $G$ has a $BN$-pair with characteristic $p$. Then $$P_{1,p,G} = {\rm St} \otimes {\rm St},$$  where ${\rm St}$ denotes the Steinberg module for $RG$.}

\medskip
\noindent {\bf Proof:} As noted in the preceding discussion, the Steinberg (genuine) module $St$ for $RG$ may be expressed in the Green ring for $RG$ as 
$${\rm St} = \sum_{J \subseteq I} (-1)^{|J|} {\rm Ind}_{P_{J}}^{G}(R),$$ where $R$ denotes the trivial module, and $I$ is the set of generating reflections for the Weyl group of $G$. Also, we have $${\rm St} = (-1)^{|I|}{\rm St}_{p}(G),$$ (according to the sign notation we have adopted).

\medskip
Whenever $P_{J}$ is a parabolic subgroup of $G$, we saw in the discussion preceding the statement of the Theorem that ${\rm Res}^{G}_{P_{J}}({\rm St})$ is the projective cover as $RP_{J}$-module of ${\rm St}_{P_{J}/U_{J}}.$

\medskip
Hence, by 11.3 above, together with the discussion preceding the statement of the Theorem, we see that in the Green ring of $RG$, we have 
$$P_{1,p,G} = \sum_{J \subseteq I}  (-1)^{|J|}{\rm Ind}_{P_{J}}^{G}(P[{\rm St}_{P_{J}/U_{J}}]) $$
$$ = {\rm St} \otimes \left[ \sum_{J \subseteq I}  (-1)^{|J|}{\rm Ind}_{P_{J}}^{G}(R) \right]$$ 
$$ = {\rm St} \otimes {\rm St} (= {\rm St}_{p}(G) \otimes {\rm St}_{p}(G)),$$ as required.

\medskip
\noindent {\bf Remark 11.7:}  This provides an alternative proof of the known fact that when $G$ is a finite group with a characteristic $p$-type $BN$-pair,   the number of $p$-elements of $C_{G}(y)$ is $|C_{G}(y)|_{p}^{2}$ for each $p$-regular $y \in G.$

\medskip
In fact, this approach gives some insight into the (well-known) fact that the value taken by the character $\chi_{I}$ afforded by the Steinberg module ${\rm St}$ at $y$  is $\pm |C_{G}(y)|_{p}$ for each $p$-regular $y \in G$, when $G$ has a $BN$-pair of characteristic $p$-type. We briefly outline the argument. If the $p$-regular element $y$  is such that $p| |C_{G}(y)|,$ then $y$ certainly normalizes a non-trivial $p$-subgroup of $G$, and then some conjugate of $y$ lies in a (proper) parabolic subgroup $P_{J}$, so we might as well suppose that $y \in P_{J}$. Now ${\rm Res}^{G}_{P_{J}}({\rm St})$ is the projective cover of ${\rm St}_{P_{J}/U_{J}}$ as $RP_{J}$-module. By induction on $|I|$, we may suppose that the Steinberg character of $P_{J}/U_{J}$ takes value $\pm |C_{P_{J}/U_{J}})(yU_{J})|_{p}$ at 
$yU_{J}$. Then the result of W.F. Reynolds quoted at the end of Section 3 (and also used in Section 7) tells us that the value of the character afforded by ${\rm St}$ at $y$ is $$|C_{U_{J}}(y)| \times \pm  |C_{P_{J}/U_{J}})(yU_{J})|_{p} = \pm |C_{P_{J}}(y)|_{p},$$ (using standard results on coprime action). This argument in fact tells us that we have $\chi_{I}(y) = \pm p^{r_{y}}$ for some non-negative integer $r_{y}$, and that $p^{r_{y}}| |C_{G}(y)|.$
On the other hand, the fact that $\chi_{I}$ is an irreducible character vanishing on $p$-singular elements of $G$ tells us that $\frac{\chi_{I}(y)}{|C_{G}(y)|_{p}}$ is an algebraic integer, so we must have $\chi_{I}(y)  = \pm |C_{G}(y)|_{p}.$ Strictly speaking, we have not yet dealt with the case that the $p$-regular element $y$ normalizes no non-trivial $p$-subgroup of $G.$ But in that case, we have $\chi_{I}(y) = \pm 1$ because no conjugate of $y$ lies in any (proper) parabolic subgroup $P_{J}.$

\medskip
Notice that the case $|I| = 1$ is covered by the above arguments. In the case $|I| = 1$, our assumption is that $G$ is a doubly transitive permutation group of degree $|U| + 1$ where $U$ is a Sylow $p$-subgroup of $G$, and the point stabilizer $B = N_{G}(U)$. In this case, $\chi_{I} = {\rm Ind}_{B}^{G}(1)-1$ is an irreducible character of degree $|U|$ and is afforded by a projective $RG$-module. If $y$ is a $p$-regular element of $G$ such that no conjugate of $y$ lies in $B$, then $\chi_{I}(y) = -1$. If the $p$-regular element $y$ lies in $B$, then ${\rm Res}^{G}_{B}(\chi_{I})$ is the character afforded by projective cover of the trivial $RB$-module, and we noted in Section 3 that this gives $\chi_{I}(y) = |C_{U}(y)|.$

\medskip
\section{More implications of conjectures for Brauer characters}%12

\medskip
Let $\Phi$ denote the $\ell \times \ell$ matrix whose $(i,j)$ entry is $\phi_{i}(y_{j})$, where $\{\phi_{i}: 1 \leq i \leq \ell \}$ are the irreducible Brauer characters of $G$ and $\{y_{j} : 1 \leq j \leq \ell \}$ is a set of representatives for the $p$-regular conjugacy classes of $G$. Let $C$ denote the Cartan matrix of $RG$, and let $Y$ denote the diagonal $\ell \times \ell$ matrix whose $i$-th entry on the main diagonal is $|C_{G}(y_{i})|.$ Let $\Theta$ denote the $\ell \times \ell$ matrix with $(i,j)$-entry $\theta_{i}(y_{j})$, where $\theta_{i}$ is the unique (Brauer) character of a projective indecomposable $RG$-module with $\langle \theta_{i}, \phi_{j} \rangle = \delta_{ij}$ for $1 \leq j \leq \ell.$ Note that if $\theta$ is the character of a projective $RG$-module, we may extend its Brauer character to a character of $G$ by decreeing that it takes value $0$ on all $p$-singular elements. Also, if $\phi$ is any Brauer character of $G$, we may consider $\phi \theta$ as a character of $G$ for the same reason. Finally, each Brauer character of a projective $\mathbb{F}G$-module may be regarded as a character of $G$ in the same fashion.

\medskip
By Brauer's orthogonality relations, we have  $\Theta = C\Phi $ and $$\overline{\Phi}^{t}C\Phi = Y.$$  Furthermore, we have 
$$C^{-1} = \Phi Y^{-1} \overline{\Phi}^{t}.$$

\medskip
Now let us consider the truncated conjugation character, which takes value $|C_{G}(x)|$  on $p$-regular elements $x$ and $0$ on all $p$-singular elements.
Notice that $$\Phi Y =\Phi \overline{\Phi}^{t}  C \Phi = \Phi \overline{\Phi}^{t} \Theta .$$ On the other hand, we have seen that the truncated conjugation character really is the character of a projective $RG$-module, and may be written as $$\sum_{i=1}^{\ell}\overline{\phi_{i}} \theta_{i}.$$ In other words,
we have $\Phi Y = M\Theta$  for some non-negative integer matrix $M$. Since $\Theta$ is an invertible matrix, we may conclude that 
$$M = \Phi\overline{\Phi}^{t}$$ and that $\Phi\overline{\Phi}^{t}$ is a matrix with non-negative integer entries.

\medskip
Suppose that now Conjecture B holds for $C_{G}(y)$ for all $p$-regular $y \in G$, so that $\Psi_{1,p,C_{G}(y)}$ is a character of $C_{G}(y)$ afforded by a projective $RC_{G}(y)$-module for each $p$-regular element $y$ of $G$. Recall that  the truncated conjugation character for $G$ is equal to $$\sum_{y \in G_{p^{\prime}}/G} {\rm Ind}_{C_{G}(y)}^{G}( \Psi_{1,p,C_{G}(y)}).$$ 

\medskip
For each $p$-regular element $y \in G$, let  $D_{y}$ be the diagonal matrix  whose $i$-th entry on the main diagonal is 
$${\rm Ind}_{C_{G}(y)}^{G}(\Psi_{1,p,C_{G}(y)})[y_{i}].$$ Note that this is the number of elements of the $p^{\prime}$-section of $y$ which commute with $y_{i}.$

\medskip
Since we are presently assuming that $${\rm Ind}_{C_{G}(y)}^{G}(\Psi_{1,p,C_{G}(y)})$$ is a character afforded by a projective $RG$-module in this case, we know that $\Phi D_{y} = M_{y}\Theta$ for some matrix $M_{y}$ with non-negative integer entries. Note also that the rank of $M_{y}$ is the number of non-zero entries on the main diagonal of $D_{y},$ which is the number of $y_{i}$ such that some conjugate of $y_{i}$ commutes with an element of the  $p^{\prime}$-section of $y$. This is simply the number of $y_{i}$ such that $C_{G}(y_{i})$ meets the conjugacy class of $y$.

\medskip
Notice that we have $$\sum_{y \in G_{p^{\prime}}/G} D_{y} = Y .$$
Now we have $$\Phi D_{y}  = M_{y} C \Phi, $$ so that $$\Phi D_{y} \Phi^{-1} = M_{y}C,$$ and 
$$\sum_{y \in G_{p^{\prime}}/G} \Phi D_{y} \Phi^{-1} =  \Phi Y \Phi^{-1} = \Phi \overline{\Phi}^{t} C\Phi \Phi^{-1} = MC.$$
Since $C$ is non-singular, we have $$\sum_{y \in G_{p^{\prime}}/G} M_{y} = M.$$

\medskip
Hence we see that $M$ is a positive definite symmetric matrix with non-negative integer entries which may be expressed in the form
$$\sum_{y \in G_{p^{\prime}}/G} M_{y},$$ where each $M_{y}$ has non-negative integer entries.

\medskip
Let us note that $$\Phi \overline{\Phi}^{t}$$ has $(i,j)$-entry $$\sum_{k=1}^{\ell} \phi_{i}(y_{k})\overline{\phi_{j}(y_{k})},$$ which is equal to $\langle \Lambda_{c,p,G}\phi_{i}, \phi_{j} \rangle.$

\medskip
This is equal (by Frobenius reciprocity) to  $$\sum_{s = 1}^{\ell} \langle {\rm Res}^{G}_{C_{G}(y_{s})} (\phi_{i})\Psi_{1,p,C_{G}(y_{s})}, {\rm Res}^{G}_{C_{G}(y_{s})}(\phi_{j}) \rangle ,$$  which is in turn equal to 
$$\sum_{s=1}^{\ell} \langle {\rm Res}^{G}_{C_{G}(y_{s})}(\phi_{i}^{\ast}), {\rm Res}^{G}_{C_{G}(y_{s})}(\phi_{j}^{\ast})\rangle. $$

\medskip
Hence, for each $s$, the matrix $M_{y_{s}}$ has $(i,j)$-entry $$\langle {\rm Res}^{G}_{C_{G}(y_{s})}(\phi_{i}^{\ast}), {\rm Res}^{G}_{C_{G}(y_{s})}(\phi_{j}^{\ast})\rangle,$$ so that $M_{y_{s}}$ is Hermitian and positive semi-definite for each $s$. As noted earlier, the rank of the matrix $M_{y_{s}}$ is the number of $i \leq \ell$ such that the conjugacy class of $y_{i}$ meets $C_{G}(y_{s}).$

\medskip
Hence we have proved:

\medskip
\noindent {\bf Theorem 12.1:} \emph{ Let $G$ be a finite group such that for each $p$-regular element $y \in G,$ the subgroup $C_{G}(y)$ satisfies Conjecture B.
Then $M = \Phi \overline{\Phi}^{t}$ may be expressed in the form $$M = \sum_{y \in G_{p^{\prime}}/G} M_{y},$$ where each $M_{y}$ is a positive semi-definite symmetric matrix with non-negative integer entries, and whose rank is the number of $p$-regular conjugacy classes of $G$ meeting $C_{G}(y).$
Furthermore, the matrices $M_{y}C$ are mutually commuting and simultaneously diagonalizable with sum $MC$ which is similar to $Y$. The eigenvalues of $M_{y}C$ are (counting repetitions) the $\ell$ values $${\rm Ind}_{C_{G}(y)}^{G}(\Psi_{1,p,C_{G}(y)})[z]$$ as $z$ runs over a set of representatives of the $p$-regular conjugacy classes of $G$. These values may be respectively rewritten as $|S_{p^{\prime}}^{G}(y) \cap C_{G}(z)|$ in each case.}

\medskip
Let us continue in a similar spirit in the case that $G$ has a $BN$-pair in characteristic $p$. Let $\chi_{s}$ denote the Steinberg character of $G,$ and recall that in this case $\Psi_{1,p,G} = \chi_{s}\overline{\chi_{s}}.$ Let $\{ \phi_{i} : 1 \leq i \leq \ell \}$ be the set of irreducible Brauer characters of $G.$ As before, we may regard
each $\phi_{i}\chi_{s}$ as the character of a projective $RG$-module by considering it to have value $0$ on all $p$-singular elements. Then we may obtain a version of a theorem of G. Lusztig [11]:

\medskip
\noindent {\bf Theorem 12.2:} \emph{ Let $G$ be as above. Then $$\{\chi_{s}\phi_{i}: 1 \leq i \leq \ell \}$$ is a $\mathbb{Z}$-basis, consisting of characters afforded by projective $RG$-modules, for the $\mathbb{Z}$-module of generalized characters of $G$ vanishing on all $p$-singular elements.}

\medskip
\noindent {\bf Proof:} Recall that $\Psi_{1,p,G} = \chi_{s}\overline{\chi_{s}}$ in this case. Now for each $i,j,$ we have 
$$\langle \phi_{i}^{\ast}, \phi_{j}^{\ast} \rangle = \langle \Psi_{1,p,G}\phi_{i},\phi_{j} \rangle = \langle \phi_{i}\chi_{s}, \phi_{j}\chi_{s} \rangle .$$

\medskip
Let $M_{1}$ be the matrix defined as in the proof of the previous theorem, so that $\Phi D_{1} = M_{1}\Theta$ and $$M_{1}C = \Phi D_{1} \Phi^{-1}.$$
Now $${\rm det} (D_{1}) = \prod_{y \in G_{p^{\prime}}/G} |C_{G}(y)|_{p}^{2} = {\rm det}(C)^{2}.$$  Hence ${\rm det}(M_{1}) = {\rm det}(C).$

\medskip
For $1 \leq i \leq \ell,$ we may write $$\chi_{s}\phi_{i} = \sum_{j = 1}^{\ell} a_{ij}\theta_{j}$$ for certain uniquely determined integers $a_{ij}$.
Let $A$ denote the $\ell \times \ell$ matrix $[a_{ij}].$ Now notice that we have $$A^{t}CA = [\langle \phi_{i}^{\ast}, \phi_{j}^{\ast} \rangle] = M_{1},$$ as in the proof of the previous theorem. Thus $A$ is unimodular and the result follows.

\section{When $p = 2$ and the Sylow $2$-subgroup of $G$ is a Klein $4$-group}%13

\medskip
In this case, $\Psi_{1,2,G}(x)$ is the number of square roots of $x$ in $G$ for each $x \in G,$  so that 
$$\Psi_{1,2,G} = \sum_{\mu \in {\rm Irr}(G)} \nu(\mu) \mu,$$ where $\nu$ is the usual Frobenius-Schur indicator. Furthermore, since $C_{G}(t)$ has a normal $2$-complement  for each involution $t \in G$ under current hypotheses, we see by Theorem 7.3 that $\Psi_{1,2,G}$ is a character of $G$. Hence $\nu(\mu) = 1$ for each real-valued irreducible character $\mu$ of $G$. 

\medskip
We now aim to prove that $\Psi_{1,2,G}$ is afforded by a projective $RG$-module. By the results of Section 7, it suffices to prove that $\Psi_{1,2,G}^{(b)}$ is afforded by a projective $RG$-module, where $b$ is the principal block of $RG$. We have\\ $[N_{G}(S):C_{G}(S)] \in  \{1, 3\},$ where $S \in {\rm Syl}_{2}(G).$ If $N_{G}(S) = C_{G}(S),$ then $G$ has a normal $2$-complement, and by Theorem 3.2, $\Psi_{1,2,G}$ is a character afforded by a projective $RG$-module. Suppose then that $[N_{G}(S):C_{G}(S)] = 3.$

\medskip
In that case, $G$ has a single conjugacy class of involutions, say with representative $t$. Furthermore, the principal $2$-block $b$ of $G$ contains four irreducible characters, all of odd degree, and we have $\ell(b) = 3.$
Let $\theta_{1}$ denote the character of the projective cover of the trivial module. Then $\theta_{1}$ is real-valued, and all its irreducible constituents lie in the block $b$.

\medskip
Let $1, \chi_{2},\chi_{3}, \chi_{4}$ be the set of irreducible characters in $b$. Then there are signs $1, \epsilon_{2},\epsilon_{3}, \epsilon_{4}$ such that $$ \beta = 1 + \epsilon_{2}\chi_{2}+\epsilon_{3}\chi_{3} + \epsilon_{4}\chi_{4}$$ vanishes on all $2$-regular elements, and $\beta$ takes value $4$ on all $2$-singular elements (note that $\epsilon_{i} = \chi_{i}(t)$ for $2 \leq i \leq 4).$ Note also that $tx$ has no square root in $G$ whenever $x$ is an element of odd order in $C_{G}(t),$ so that  $1 + \sum_{i=2}^{4}\epsilon_{i} \nu(\chi_{i}) =0,$ by Brauer's Second and Third Main Theorems.

\medskip
Suppose first that $\chi_{2}$ is not real-valued. Then $$\sum_{g \in G} \chi_{2}(g) = 0.$$ Also, we know that $$\sum_{g \in G} \chi_{2}(g^{2}) = 0$$ since $\chi_{2}$ has Frobenius-Schur indicator $0$. This leads to $$\frac{|G|}{4} \left( \langle {\rm Res}^{G}_{O_{2^{\prime}}(C_{G}(t))}(\chi_{2}), 1 \rangle \right) = \epsilon_{2}\frac{|G|}{4},$$ which forces $\epsilon_{2}$ to be non-negative, and hence $\epsilon_{2} =1$. But in that case, we may also assume that $\chi_{3} = \overline{\chi_{2}}$ and so also $\epsilon_{3} = 1.$ This now forces $\epsilon_{4} = -1$ and $\nu(\chi_{4}) = 1.$

\medskip
In this case, then, we have $\Psi_{1,2,G}^{(b)} = 1 + \chi_{4}.$

\medskip
Now suppose that $\chi_{i}$ is real-valued for $2 \leq i \leq 4.$ Then $\chi_{i}$  has Frobenius-Schur indicator $1$ for $1 \leq i \leq 4$, and we have $$\Psi_{1,2,G}^{(b)} = 1 + \chi_{2} + \chi_{3} + \chi_{4}.$$

\medskip
By results of Erdmann [5,6] and Landrock [10] (and earlier results of Brauer) we know that all decomposition numbers for $b$ are $0$ or $1$.
Also, up to ordering of simple $b$-modules, there are only two possibilities for the Cartan matrix of $b$. These are 
$$\left( \begin{array}{clcr} 4 &2 &2\\2&2 &1\\2&1&2 \end{array} \right)$$ and $$\left( \begin{array}{clcr} 2 &1 &1\\1&2 &1\\1&1&2 \end{array} \right).$$
In the former case, where $b$ is Morita equivalent to the principal $2$-block of $A_{5},$ the character of the projective cover of the trivial module has norm-squared 4. In the second case where $b$ is Morita equivalent to the principal $2$-block of $A_{4},$ the character of the projective cover of the trivial module  has norm squared $2.$

\medskip
In all cases, we may conclude that $\Psi_{1,2,G}^{(b)} = \theta_{1},$ the character afforded by the projective cover of the trivial module. For if $\langle \theta_{1},\theta_{1} \rangle = 4,$ and $$\Psi_{1,2,G}^{(b)} = 1 + \chi_{4},$$ then we are forced to conclude that $\chi_{2} + \chi_{3}$ vanishes on $2$-singular elements when $\chi_{2}$ is not real valued. However, $\chi_{2}$ and $\chi_{3}$ agree and do not vanish on $2$-singular elements in that case, a contradiction.

\medskip
On the other hand, in the case that $\chi_{2}, \chi_{3},\chi_{4}$ are all real-valued and $$\theta_{1}  = 1+\chi_{2},$$ (without loss of generality), then we may label so that  $\epsilon_{2} = -1 = \epsilon_{4} = - \epsilon_{3},$ since $1 + \epsilon_{2} + \epsilon_{3} + \epsilon_{4} = 0.$ But the only characters of projective indecomposable $RG$-modules in $b$ which are consistent with the Cartan matrix for $b$ in this case are $1+ \chi_{2}, \chi_{2} + \chi_{3}$ and $\chi_{2}+ \chi_{4},$ and these do not all vanish at $t$. Hence we have proved:

\medskip
\noindent {\bf Theorem 13.1:} \emph{ Suppose that $G$ is a finite group whose  Sylow $2$-subgroup is a Klein $4$-group. Then $\Psi_{1,2,G}$ is a character afforded by a projective $RG$-module.}

\medskip
In particular, we have:

\medskip
\noindent {\bf Corollary 13.2:} \emph{ Suppose that $G \cong {\rm PSL}(2,q)$ with $q \not \equiv \pm 1$ (mod $8$). Then Conjecture B holds for $G$ 
(for each prime divisor $p$ of $|G|$).}

\medskip
\noindent {\bf Proof:} If $q =2,$ then $G$ has a normal $2$-complement and a normal Sylow $3$-subgroup. If $q = 3,$ then $G$ has a normal  $3$-complement and a normal Sylow $2$-subgroup. Theorem 3.2 may be applied in these cases.

\medskip
If $q >3,$ then $G$ is a finite simple group. Let $p$ be a prime divisor of $|G|$. If $p$ divides $q,$ then $G$ has a characteristic $p$ type BN-pair and $\Psi_{1,2,p}$ is the square of the Steinberg character. Hence we may suppose that $q > 3$ and $p$ does not divide $q$. If $p$ is odd, then $G$ has a cyclic Sylow $p$-subgroup, so that $\Psi_{1,p,G}$ is a character afforded by a projective $RG$-module. If $p =2,$ then $G$ has a Klein $4$ Sylow $2$-subgroup and the result follows by 13.1.

\section{Conjecture $B$ for $G = {\rm PSL}(2,q)$ when ${\rm PSL}(2,q)$ has order divisible by $8$}%14

\medskip
Throughout this section, we suppose that $G \cong {\rm PSL}(2,q)$ with $q$ odd and that $G$ has a dihedral Sylow $2$-subgroup with at least $8$ elements. Notice that $G$ has one conjugacy class of involutions in this case, say with representative $t$. Our aim is to prove that $\Psi_{1,2,G}$ is a character of $G$ which may be  afforded by a projective $RG$-module. This will complete the proof that Conjecture B holds for $G$ for all primes when $G \cong {\rm PSL}(2,q)$. By the results of Sections 7,9 and 13, it only remains to deal with the case $p=2$ when $G$ has a dihedral Sylow $2$-subgroup of order at least $8$, and we need to prove that whenever $\phi_{i}$ is an absolutely irreducible Brauer character in the principal $2$-block of $RG$, then $$\langle \Psi_{1,2,G}, \phi \rangle \geq 0.$$

\medskip
There are two cases to consider: 

\medskip
\noindent{\bf Case 1:} $q \equiv 1$ (mod $8$). In this case, $C_{G}(t)$ is a dihedral group with $q-1$ elements, and has a unique cyclic subgroup $H$ of index $2$. Write $H = K \times L$ where $K$ has odd order and $L$ is a $2$-group. Then $K = O_{2^{\prime}}(C_{G}(x))$ for each non-identity element $x \in L.$

\medskip
The character table of $G$ was known to Frobenius. Combining this with the results of Brauer, Erdmann and Landrock mentioned in the previous section, we have the following  information: $G$ has exactly four irreducible characters of odd degree (respectively of degree $1,q, \frac{q+1}{2}, \frac{q+1}{2}$), and these all lie in the principal $2$-block $b$ of $G$ (in general, as noted by R. Brauer, any $2$-block of defect greater than one contains at least $4$ irreducible characters of odd degree). 

\medskip
$G$ has $\frac{q-1}{4}$ irreducible characters of degree $q-1$, which all lie in $2$-blocks of defect zero of $G$.

\medskip
$G$ has $\frac{q-5}{4}$ irreducible characters of degree $q+1$, and $\frac{|L|}{2} -1$ of  these are (all) the height one irreducible characters in the principal $2$-block of $G$.

\medskip
By the results of Section VII of (Burkhardt, [3]) we know in this case that the Brauer character (for the prime $2$) of a (complex) irreducible character $\chi$  of degree $\frac{q+1}{2}$ decomposes on reduction (mod $J(R)$ ) as $1 + \phi$ where $\phi$ is a Brauer irreducible character of degree $\frac{q-1}{2}.$ We claim that  $$\langle \Psi_{1,2,G}, \phi \rangle \geq 0.$$ By Theorem 7.3, we know that $$\langle \Psi_{1,2,G}, \chi \rangle \geq 0,$$ so we only need to exclude the possibility that $\langle \Psi_{1,2,G}, \chi \rangle = 0,$ since we know that $\langle \Psi_{1,2,G}, 1 \rangle = 1.$

\medskip
The argument at the end of the proof of 7.3 shows that the only way we can obtain $\langle \Psi_{1,2,G}, \chi \rangle = 0$ is if we have 
$$\langle {\rm Res}^{G}_{O_{2^{\prime}}(C_{G}(x))}(\chi),  1 \rangle = \chi(x)$$ for each non-identity $2$-element $x \in G.$
But in this case, $O_{2^{\prime}}(C_{G}(t))$ is cyclic of order at most $\frac{q-1}{8}$ for an involution $t \in G$, and $\chi$ takes constant value $1$ on 
non-identity elements of $O_{2^{\prime}}(C_{G}(t))$, (because the Borel subgroup $B$ of ${\rm PSL}(2,q)$ is a Frobenius group of order $\frac{q(q-1)}{2}$ and, (using unimodularity),  ${\rm Res}^{G}_{B}(\chi)$ decomposes as $1 + \mu,$ where $\mu$ is irreducible of degree $\frac{q-1}{2}$ and $\mu$ vanishes on all non-identity elements of the cyclic Frobenius complement). Hence
 $$\langle {\rm Res}^{G}_{O_{2^{\prime}}(C_{G}(t))}(\chi),  1 \rangle \geq	 \frac{8(q+1)}{2(q-1)} > 4,$$ so we certainly have $$\langle {\rm Res}^{G}_{O_{2^{\prime}}(C_{G}(t))}(\chi),  1 \rangle  > 1 = \chi(t)$$ and hence $\langle \Psi_{1,2,G}, \chi \rangle \geq  1,$ as required to prove the claim.

\medskip
Since it was known to Brauer that there are three absolutely irreducible Brauer characters in the principal $2$-block of $G$ in this case, and from [3], we know that these are the trivial Brauer character and two algebraically conjugate irreducible Brauer characters of degree $\frac{q-1}{2},$ we may conclude in this case that $\Psi_{1,2,G}$ has non-negative inner product with each Brauer irreducible character of $G$, so that $\Psi_{1,2,G}$ may be afforded by a projective $RG$-module.

\medskip
\noindent {\bf Case 2 :} $q \equiv -1$ (mod $8$).

\medskip
In this case, we know from [3] again that there are three absolutely irreducible Brauer characters in the principal $2$-block of $G$, and they again have respective degrees $1,\frac{q-1}{2},\frac{q-1}{2}.$ However, in this case, these are all obtained by reduction (mod $J(R)$) of complex irreducible characters of $G$ (this time, the odd degree irreducible characters in the principal $2$-block have respective degrees $1,q, \frac{q-1}{2},\frac{q-1}{2}$).
By Theorem 7.3, it follows that $\langle \Psi_{1,2,G}, \phi \rangle \geq 0 $ for each Brauer irreducible character $\phi$, so that we have completed the proof that $\Psi_{1,2,G}$ may be afforded by a projective $RG$-module.

\medskip
\noindent {\bf Theorem 14.1:}  \emph{Suppose that $G \cong {\rm PSL}(2,q)$ where $q$ is a prime power. Then  for each prime $p$, the generalized character $\Psi_{1,p,G}$ is a character which may be afforded by a projective $RG$-module.}

\medskip
In fact, using Theorem 3.2, we may deduce:

\medskip
\noindent {\bf Corollary 14.2:} \emph{Suppose that $G \cong {\rm SL}(2,q)$ where $q$ is a prime power. Then  for each prime $p$, the generalized character $\Psi_{1,p,G}$ is a character which may be afforded by a projective $RG$-module.}

\medskip
\noindent {\bf Proof:} If $p$ is an odd prime divisor of $|{\rm SL}(2,q)|,$ then $\Psi_{1,p,G}$ contains the central involution of ${\rm SL}(2,q)$ in its kernel,
and the result follows by 14.1. If $p =2,$ then by Theorem 3.2 we may deduce that $\Psi_{1,2,G}$ is afforded by the projective cover of $M$ as $RG$-module, where $\Psi_{1,2,G/Z(G)}$ is afforded by the projective $RG/Z(G)$-module $M$.

\section{The truncated conjugation module as a direct summand of copies of the regular module}%15

\medskip
\noindent {\bf Theorem 15.1 :} \emph{ The truncated conjugation module $T_{c,p,G}$ is a direct summand of the direct sum of $\ell$ copies of the regular $RG$-module, and is not a direct summand of the direct sum of $r$ copies of the regular module for any $r < \ell.$}

\medskip
\noindent {\bf Proof:} We know that the projective cover of the trivial module occurs with multiplicity $\ell$ as a summand of $T_{c,p,G}$, while it only occurs with multiplicity one as a direct summand of the regular module $RG$. Hence the second assertion is clear. To prove the first assertion, it suffices to pove that for $1 \leq i \leq \ell,$ the character $\theta_{i}$ occurs with multiplicity at most $\ell \phi_{i}(1)$ as a summand of $\Lambda_{c,p,G}$. But this multiplicity is $\sum_{j = 1}^{\ell} \phi_{i}(y_{j})$, which is a non-negative integer less than or equal to $\ell \phi_{i}(1).$

\medskip
\noindent {\bf Corollary 15.2:} \emph{ Whenever $P_{1,p,G}$ is a projective $RG$-module, then it is a direct summand of the regular $RG$-module.}

\medskip
\noindent {\bf Proof:} It suffices to prove that for $i \leq \ell,$ the multiplicity of $\theta_{i}$ in $\Psi_{1,p,G}$ is at most $\phi_{i}(1).$
 This multiplicity is $$\sum_{j=1}^{\ell} \frac{\phi_{i}(y_{j}) |G_{p} \cap C_{G}(y_{j})|}{|C_{G}(y_{j})|}.$$ Since the multiplicity of the character $\theta_{1}$ of the projective cover of the trivial module in $\Psi_{1,p,G}$ is one, we see easily that the multiplicity of $\theta_{i}$ is at most $\phi_{i}(1)$, and that equality occurs if and only if $O^{p}(G) \leq {\rm ker} \phi_{i},$ in other words, equality can only occur when $i = 1$ and $\phi_{i}$ is the trivial Brauer character. Hence we have proved that $\theta_{i}$ occurs with multiplicity at most $\phi_{i}(1) - 1$ in $\Psi_{1,p,G}$ for $i > 1.$

\medskip
\noindent {\bf Remark 15.3:} In particular, whether or not $P_{1,p,G}$ is a projective $RG$-module, the only character afforded by the projective cover (as $RG$-module) of a one-dimensional simple $\mathbb{F}G$-module which occurs with non-zero multiplicity in $\Psi_{1,p,G}$ is the character afforded by the projective cover of the trivial module.

\medskip
\noindent {\bf Acknowledgement:} We are grateful to G. Malle for pointing out some typographical errors and other inaccuracies in an earlier version of these notes, as well as providing the reference [3].

\medskip
\begin{center} {\bf References} \end{center}
\noindent [1] Boltje, R., \emph{A canonical Brauer induction formula}, Ast\'erisque, {\bf 181-182}, (1990), 31-59.

\medskip
\noindent [2]  Brauer, R. and Feit, W., \emph{On the number of irreducible characters in a given block}, Proceedings of the National Academy of Sciences,
{\bf 45},3,(1959), 361-365.

\medskip
\noindent [3] Burkhardt, R., \emph{Die Zerlegungsmatrizen der Gruppen ${\rm PSL}(2,p^{f})$}, Journal of Algebra, {\bf 40}, (1976), 75-96.

\medskip
\noindent [4] Curtis, C.W., \emph{ The Steinberg character of a finite group with a $(B,N)$-pair}, Journal of Algebra, {\bf 4}, 3,(1966), 433-441.

\medskip
\noindent [5] Erdmann, K., \emph{Principal blocks of groups with dihedral Sylow $2$-subgroups}, Communications in Algebra,{\bf 5},(1977), 665-694.

\medskip
\noindent [6] Erdmann, K., \emph{Blocks whose defect groups are Klein $4$-groups}, Journal of Algebra, {\bf 59},(1979),452-465.

\medskip
\noindent [7] Fong, P., \emph{Solvable groups and modular representation theory}, Transactions of the American Mathematical Society, {\bf 103},3,(1962),484-494.

\medskip
\noindent [8] Green, J.A., \emph{The characters of the finite general linear group}, Transactions of the American Mathematical Society, {\bf 80},2,(1955),403-447.

\medskip
\noindent [9] Isaacs, I.M., \emph{Characters of $\pi$-separable groups}, Journal of Algebra, {\bf 86},(1984),98-128.

\medskip
\noindent [10] Landrock, P., \emph{The principal block of finite groups with dihedral Sylow $2$-subgroups}, Journal of Algebra,{\bf 39},2,(1970),410-428.

\medskip
\noindent [11] Lusztig, G., \emph{Divisibility of projective modules of finite Chevalley groups by the Steinberg module}, Bulletin of the London Mathematical Society, {\bf 8}, 2, (1976), 130-134. 

\medskip
\noindent [12] Murray, J., \emph{Projective modules and involutions}, Journal of Algebra,{\bf 299}, 2,(1006), 616-622.

\medskip
\noindent [13] Reynolds, W.F., \emph{Blocks and normal subgroups of finite groups}, Nagaoya Mathematical Journal, {\bf 22}, (1963), 15-32.

\medskip
\noindent [14] Robinson, G.R., \emph{On projective summands of induced modules}, Journal of Algebra, {\bf 122}, (1989), 106-111.

\medskip
\noindent [15] Scharf, T., \emph{Die Wurzlanzahlfunktion in symmetrichen Gruppen}, Journal of Algebra, {\bf 131},(1991),446-457.

\medskip
\noindent [16] Solomon, L., \emph{On the sum of the elements in the character table of a finite group}, Proceedings of the American Mathematical Society, {\bf 12},(1961),961-963.

\medskip
\noindent [17] Webb, P.J., \emph{Subgroup complexes}, {\bf in}: The Arcata Conference on Representations of Finite Groups (\emph{ed. P. Fong}), AMS Proceedings of Symposia in Pure Mathematics, {\bf 47}, (1987).

\end{document}